\documentclass[11pt]{article}
\usepackage{a4wide}
\usepackage{amsmath,amsthm,amssymb,amsfonts,esint}
\usepackage{graphicx}
\usepackage[usenames,dvipsnames]{color}
\usepackage[colorlinks=true, pdfstartview=FitV, linkcolor=blue, citecolor=blue, urlcolor=blue]{hyperref}
{\newtheorem{thm}{Theorem}[section]}
{}
{}
{\newtheorem{lemme}[thm]{Lemma}}
{}
{\newtheorem{rem}{Remark}[section]}

\newcommand{\na}{\nabla}
\newcommand{\eps}{\epsilon}
\newcommand{\pa}{\partial}
\newcommand{\N}{{\mathbb N}}
\newcommand{\R}{{\mathbb R}}
\newcommand{\T}{{\mathbb T}}
\newcommand{\mP}{{\mathbb P}}
\newcommand{\mQ}{{\mathbb Q}}

\renewcommand{\div}{{\mathrm{div} \,}}
\newcommand{\curl}{{\mathrm{curl} \,}}

\begin{document}

\title{Diffusion-free boundary conditions \\ for the Navier-Stokes equations}
\author{Emmanuel Dormy\thanks{emmanuel.dormy@ens.fr, Ecole Normale Supérieure, CNRS, DMA, 75005 Paris, France}, David Gerard-Varet\thanks{david.gerard-varet@imj-prg.fr, Université Paris Cité and Sorbonne Université, CNRS, IMJ-PRG, F-75013 Paris, France}}
\maketitle

\begin{abstract}
    We provide a mathematical analysis of the `diffusion-free' boundary conditions recently introduced  by Lin and Kerswell \cite{Kerswell24b} for the numerical treatment of inertial waves in a fluid contained in a rotating sphere. We consider here the full setting of the nonlinear Navier-Stokes equation in a general bounded domain $\Omega$ of $\R^d$, $d=2$ or $3$. We show that diffusion-free boundary conditions
    \begin{eqnarray*}
&         \Delta u \cdot \tau \vert_{\pa \Omega} = 0, & \quad u \cdot n\vert_{\pa \Omega} = 0 \quad \text{ when } d=2, \\
 &        \Delta u \times n\vert_{\pa \Omega} = 0, & \quad u \cdot n\vert_{\pa \Omega} = 0 \quad \text{ when } d=3, 
    \end{eqnarray*}
    allow for a satisfactory well-posedness theory of the full Navier-Stokes equations (global in time for $d=2$, local for $d=3$).  Moreover, we perform a boundary layer analysis in the limit of vanishing viscosity $\nu \rightarrow 0$. We establish that the amplitude of the boundary layer flow is in this case of order $\nu$, i.e. much lower than in the case of standard Dirichlet or even stress-free conditions. This confirms analytically that this choice of boundary conditions may be used to reduce diffusive effects in numerical studies relying on the Navier-Stokes equation to approach nearly inviscid solutions. 
\end{abstract}

\section*{Introduction}

The motion of fluids is usually described, under the approximation of incompressibility, using the 
Navier-Stokes equations 
\begin{equation} \label{NS}
    \begin{aligned}
        \pa_t u + u \cdot \na u + \na p - \nu \Delta u & = 0 \quad \text{in} \quad \Omega, \\
        \div u & = 0  \quad \text{in}\quad \Omega \, .
    \end{aligned}
\end{equation} 
In the case of large scale geophysical or astrophysical flows, this system is often completed by the Coriolis force. 
Most of the flows in nature are very large scale and therefore only weakly influenced by viscosity.
Except for particular situations, such as boundary layers detachment, the large scale flow away from boundaries is expected to be well described by the inviscid solution, i.e. by the Euler equation.
This is true of the large scale flows in the atmosphere, in the ocean, in the core of the Earth, or in the interior of stars. The numerical resolution of the Navier-Stokes equations in such a limit of vanishing viscosity causes significant challenges. Because of practical computational limitations, viscous effects are often over-estimated in numerical models. In other words, numerical simulations are often, if not always, performed with over estimated values of the viscosity. The main effect on the solution appears at the boundaries via enhanced boundary layers, which causes significant deviation of the numerical work from the physical picture. Another case in which the formation of boundary layers is to be avoided in numerical simulations is that of artificial boundaries introduced to close unbounded domains. 

Using for instance classical `no slip' (NS) boundary conditions (i.e. homogeneous Dirichlet boundary conditions) at the boundary of the domain,
\begin{equation} \label{Noslip}
    u \vert_{\pa \Omega} = 0
\end{equation}
yields in particular too strong boundary layers that significantly dissipate the energy of the main flow. 
From a mathematical point of view, these boundary layers are characterised by a ${\cal O}(1)$ amplitude correction, as such, they are very unstable. This is in part the reason why the convergence of smooth Navier-Stokes solutions to Euler solutions still is an open problem.

Since the real parameters of the natural flows are usually out of  reach of simulations, numericists have developed strategies to approach the relevant limit of vanishing viscosity with finite computational resources. 
 The most classical approach to reduce viscous effects at the boundary and thus approach more efficiently the limit of vanishing viscous effects is known as using `stress-free' (SF), also sometimes called `free-slip', boundary conditions. The impermeability condition $u\cdot n =0$ is then preserved, but the solid is assumed to exert no tangential stress on the fluid. In dimension $d=3$, this reads
\begin{equation}
    u \cdot n\vert_{\pa \Omega} = 0 \, , \quad {\rm and} \quad D(u)n \times n  \vert_{\pa \Omega} = 0 \, , 
    \label{BC_SF1}
\end{equation}
where $D(u) = \frac{\nabla u + \nabla u^t}{2}$ refers to the symmetric part of $\nabla u$, while $n$ is the unit normal vector {\em pointing  outward $\Omega$.}
This choice of boundary conditions significantly reduces viscous damping near the solid walls, thus allowing for a faster convergence to the limiting solution of vanishing viscosity. Internal stresses however exist within the flow, and with this choice of boundary condition, these need to vanish at the boundary. A boundary layer thus still forms in the flow, albeit much weaker that  in the `no slip' case. The amplitude of the boundary layer correction is now ${\cal O}(\sqrt{\nu}) \, $, compared to ${\cal O}(1)$ with NS conditions. See for instance the analysis in \cite{IfSu2011,GiKe2012}, where the more general case of Navier type boundary conditions, involving a slip length $\lambda$, is considered
\begin{equation}
   u \cdot n\vert_{\pa \Omega} = 0 \, , \quad {\rm and} \quad -\lambda \, D(u)n \times n  \vert_{\pa \Omega} = u \times n  \vert_{\pa \Omega} \, .  
\end{equation}

With the choice of SF boundary conditions, the solid does not exert any tangential stress on the fluid. As such, if the domain is sufficiently simple, the fluid conserves its linear and angular momentum. For example in a pipe, the initial linear momentum will be conserved, and in an axisymmetric domain such as a sphere, the initial angular momentum will be conserved. 

The SF condition is equivalent to 
\begin{equation}
    u \cdot n\vert_{\pa \Omega} = 0 \, ,\quad {\rm and} \quad \pa_n u \times n  \vert_{\pa \Omega} - \kappa (u \times n)\vert_{\pa \Omega} = 0 \, ,
    \label{BC_SF2}
\end{equation}
where $\kappa = \kappa(x) =  Dn(x)$ is the Weingarten endomorphism: see \cite[Lemma 1, p. 233]{BDGV2007}. As $\pa_n u = D(u)n + \frac{1}{2} \omega \times n$, where $\omega = \nabla \times u$ is the vorticity, we can re-express the boundary condition as   
\begin{equation} 
    u \cdot n\vert_{\pa \Omega} = 0 \, ,\quad {\rm and} \quad \omega \times n  \vert_{\pa \Omega}  - \kappa (u \times n)\vert_{\pa \Omega} = 0 \, .
    \label{BC_SF3}
\end{equation}
An alternative modified boundary condition was introduced by Lions \cite[pp. 87-98]{LionsJL}. It also reduces viscous effects at the boundary and takes the form 
\begin{equation}
    u \cdot n\vert_{\pa \Omega} = 0 \, , \quad {\rm and} \quad  \omega \times n  \vert_{\pa \Omega} = 0 \, .
    \label{BC_Lions}
\end{equation}
By comparison with \eqref{BC_SF3}, we see that Lions' conditions reduces to stress-free conditions for flat boundaries, but differs otherwise. For further discussion on these boundary conditions and related ones, we refer to  \cite{Kel2006,XiXi2007}.

Recently, an alternative set of boundary conditions, called `diffusion-free' has been numerically tested on several problems
\cite{Kerswell24a, Kerswell24b} from geophysics or non-newtonian flows. In the context of the classical Navier-Stokes equation, this amounts to maintaining impermeability but simply dropping viscous effects at the boundary on the tangential components, i.e.
\begin{subequations} \label{BC}
\begin{align}
&         \Delta u \cdot \tau \vert_{\pa \Omega} = 0, \quad \: \: u \cdot n\vert_{\pa \Omega} = 0 \quad \text{ when } d=2, \label{BC2D} \\ 
 &        \Delta u \times n\vert_{\pa \Omega} = 0, \quad u \cdot n\vert_{\pa \Omega} = 0 \quad \text{ when } d=3. \label{BC3D}
 \end{align}
\end{subequations} 
Note that in dimension $d=2$, $(\tau,n)$ refers to the usual direct Frenet basis. 
%
%
Our objective in this paper is to provide a mathematical basis for the relevance and efficiency of these boundary conditions. 
As detailed in the next section, we shall consider two aspects of the mathematical theory of \eqref{NS}-\eqref{BC}. First, we will study the well-posedness of \eqref{NS}-\eqref{BC}  in a bounded domain $\Omega$. In general terms, we will build strong solutions  on some time interval $(0,T)$, with $T$  arbitrary for $d=2$ (but $T = T_\nu$ small for $d=3$).  Then, in the case $\Omega = \R^2_+$, we will perform a boundary layer analysis of the system.  We will first show that for any given Euler velocity field $u^0$ in $\Omega$,  there exists for any $N$ a $O(\nu^N)$ approximate solution $u^\nu_{app}$ of \eqref{NS}-\eqref{BC}, globally in time, behaving asymptotically as 
\begin{equation}
u^\nu_{app} \approx u^0 + \nu u^{bl}(t,x,y/\sqrt{\nu}) + o(\nu) \, .
\end{equation}
We will then prove that on any time interval $(0,T)$, any strong solution $u^\nu$ of \eqref{NS}-\eqref{BC} that starts close to $u^\nu_{app}$ remains close. 

Surprisingly, dropping the highest order term in the governing equations \eqref{NS} at the boundary for the tangential flow \eqref{BC} provides a sensible set of boundary conditions and thus offers a strategy to numerically approach the limit of vanishing viscosity in a flow governed by the Navier-Stokes equations. Further discussion of advantages and drawbacks of the diffusion-free boundary conditions will be provided in the final section.

\section{Statement of the results} 

\subsection{Well-posedness theory}
The well-posedness theory of system \eqref{NS}-\eqref{BC} will be addressed  in Sections \ref{secWP2d} and \ref{secWP3d}. We will focus on the case of a smooth simply connected bounded domain $\Omega \subset \R^d$, $d=2$ or $d=3$. Extension to a half-plane, resp. to non-simply connected 2d domains, will be discussed in Paragraph \ref{halfplane}, resp. Paragraph \ref{nonconnected}. In both dimensions 2 and 3, we find it convenient to work with the vorticity formulation, in contrast to \cite{Kerswell24b} which relies on a poloidal-toroidal decomposition of the velocity field. 

\paragraph{2D setting.} The two-dimensional setting is more favorable, with a global in time well-posedness result. 
The vorticity $\omega = \curl u = \pa_x u_y - \pa_y u_x$ obeys the equation 
\begin{equation} \label{eqvorticity2d}
\pa_t \omega + u \cdot \na \omega - \nu \Delta \omega = 0 \quad \text{in } \Omega, \quad u = \mathcal{BS}\omega \, ,
\end{equation}
where $\mathcal{BS}$ denotes the Biot-Savart operator. It is useful to recall that $\mathcal{BS}f$ is the unique velocity field 
satisfying 
$$ \curl\mathcal{BS}f = f\quad \text{in } \Omega, \quad \mathcal{BS}f   \cdot n\vert_{\pa \Omega} = 0.  $$
As $\Omega$ is simply connected, $\mathcal{BS}$ is well-defined as an operator from $W^{k,p}(\Omega)$ to $W^{k+1,p}(\Omega)$ for all $1 < p < \infty$, for all $k \in \N$. Note that  $u$ can be expressed in terms of a stream function $\psi$, $u = \na^\perp \psi\equiv (-\partial _y \psi, \partial_x \psi )$, with $\Delta \psi = \omega$ in $\Omega$, $\psi\vert_{\pa \Omega}$, so that $\mathcal{BS} = \na^\perp \Delta^{-1}$, where $\Delta^{-1}$ refers to the Dirichlet Laplacian. 

A key point for the analysis of  \eqref{BC2D} is that it  can be reformulated as 
\begin{equation} \label{BCvorticity2d}
\pa_n \omega\vert_{\pa \Omega} = 0, \quad  u \cdot n\vert_{\pa \Omega} = 0 \, .
\end{equation}
This follows  from the identity: $\Delta u =  \na^\perp \omega \, ,$ so that $\Delta u \cdot \tau  = \pa_n \omega$. This Neumann condition simplifies significantly the derivation of good stability estimates, resulting in the following theorem:
\begin{thm} \label{thmWP2d}   
 Let $\Omega$ a $C^{1,1}$ simply connected bounded  domain. Let $T > 0$, $\omega_0 \in L^2(\Omega)$. There exists a unique solution $\omega$ of \eqref{eqvorticity2d}-\eqref{BCvorticity2d} such that, for all $\theta > 0$, for some $C > 0$ independent of $\nu \le 1$,
\begin{equation} \label{estim2d1}
\sup_{0 < t \le T} \|\omega(t)\|^2_{L^2} +  \nu^\theta \int_0^T \|\frac{d\omega}{dt}(s)\|^2_{(H^1)'} ds +  \nu \int_0^T \|\omega(s)\|_{H^1}^2 ds \le C \|\omega_0\|_{L^2}^2, \quad \omega\vert_{t=0} = \omega_0. 
\end{equation}
Moreover, in the case $\omega_0 \in H^1(\Omega)$, we get the following additional estimate,  for some $C_{0,\nu}$ (depending on $\|\omega_0\|_{H^1}$ and $\nu$)
\begin{equation} \label{estim2d2}
\sup_{0 < t \le T} \|\omega(t)\|^2_{H^1}  +  \int_0^T \|\frac{d\omega}{dt}(s)\|^2_{L^2} ds +  \int_0^T \|\omega(s)\|_{H^2}^2 ds \le C_{0,\nu}.
\end{equation}
\end{thm}    
\begin{rem} \label{rem_initial_boundary}
(Meaning of the initial and the boundary conditions)   It follows from \eqref{estim2d1} that $\omega$ is continuous in time with values in $(H^1)'$ so that equality $\omega\vert_{t=0} = \omega_0$ makes sense at least in $(H^{-1})'$. Regarding the Neumann condition on $\omega$, it is as usual understood in a variational way: for any $\varphi \in H^1(\Omega)$,
$$ \langle \Delta \omega, \varphi \rangle_{\langle H^{-1}, H^1\rangle} = - \int_{\Omega} \na \omega \cdot \na \varphi$$
\end{rem}
\begin{rem} \label{rem_Leray}
The boundary condition \eqref{BC2D}, rewritten as \eqref{BCvorticity2d}, complements nicely the vorticity formulation of the equation. This naturally prompts us to consider initial vorticities that are at least in $L^2(\Omega)$. Existence of Leray solutions, starting  from velocities $u_0 \in L^2(\Omega)$ only, is open: even making sense of the boundary condition at this level of regularity is unclear. 
\end{rem}
\begin{rem} \label{rem_Rousset}
    While estimate \eqref{estim2d1} is uniform in $\nu$, the second estimate \eqref{estim2d2} is not. We believe this limitation to be only  technical: as seen in Theorem \ref{thm_BL}, the velocity in the boundary layer associated with \eqref{BCvorticity2d} has  size $\sqrt{\nu}$ and amplitude $\nu$, so that uniform $H^1$ bound for the vorticity is expected. One may hope to obtain more accurate estimates through the use of higher order conormal derivatives, in the spirit of the analysis performed by Masmoudi and Rousset with Navier boundary conditions \cite{MasRou}. This is an interesting open problem. 
\end{rem}

\paragraph{3D setting.} 
In the 3D case, it is also convenient to reformulate the system in terms of the vorticity field $\omega = \curl u = \nabla \times u$. The Navier-Stokes equation yields
\begin{equation} \label{eqvorticity3d}
\pa_t \omega + u \cdot \na \omega  - \nu \Delta \omega = \omega \cdot \na u, \quad \div \omega = 0, \quad   u = \mathcal{BS} \omega \, ,
\end{equation}
with an extra stretching term on the right-hand side. Then, as $\Delta u = - \curl \omega$, the boundary conditions \eqref{BC3D} take the form 
\begin{equation} \label{BCvorticity3d}
   \curl \omega \times n\vert_{\pa \Omega} = 0, \quad  u \cdot n\vert_{\pa \Omega} = 0 \, .
\end{equation}
Note that in \eqref{eqvorticity3d}, as in the 2D setting, one can reconstruct the velocity field $u$ from $\omega$ by a Biot-Savart law \cite{Cheng}. More precisely, given a smooth simply connected domain $\Omega$, for any given $f$ satisfying the compatibility conditions
\begin{equation} \label{compatibilityBS}
     \div f = 0 \quad \text{in } \Omega, \quad \int_{\Gamma} f \cdot n = 0 \, ,
\end{equation} 
for any connected component $\Gamma$ of $\pa \Omega$, there exists a unique field $\mathcal{BS}f$ solution of
$$\curl \mathcal{BS}f = f, \quad \div \mathcal{BS}f = 0   \quad \text{in } \: \Omega, \quad \mathcal{BS}f\cdot n \vert_{\pa \Omega} = 0,$$
with the estimates 
$$\forall s, \quad \|\mathcal{BS}f\|_{H^{s+1}(\Omega)} \le C_s \|f\|_{H^s(\Omega)} \, .$$
The main difference with the 2D setting is the need for conditions \eqref{compatibilityBS}. These conditions are necessary for $f$ to be the curl of a vector field: see \cite[Theorem 1.1 and remark 1.2]{Cheng}. In the context of \eqref{eqvorticity3d}, as $\omega$ is divergence-free and the domain $\Omega$ is simply connected, $\omega$ can be written as a $\curl$  and therefore satisfies \eqref{compatibilityBS}. This allows to reconstruct the velocity as $u = \mathcal{BS} \omega$.

As we will show in section \ref{secWP3d}, combination of the first boundary condition in \eqref{BCvorticity3d} and of the enstrophy dissipation term in \eqref{eqvorticity3d} will yield a control of $\curl \omega$ in $L^2$. However,  while for $d=2$ this gives a control of $\omega$ in $H^1$, this is no longer the case for $d=3$. Showing the existence of solutions is thus more complicated, and requires more regularity on the data. Let $\mathbb{P}$ the Leray projector, that is the orthogonal projector in $L^2(\Omega)$ on the subspace $L^2_\sigma(\Omega) = \{v \in L^2(\Omega), \: \div v = 0, \: v\cdot n\vert_{\pa \Omega} = 0 \}$. We prove 
\begin{thm}  \label{thmWP3d}
 Let $\omega_0 \in H^2(\Omega)$, divergence-free. There exists $T_\nu > 0$, $C_{0,\nu} > 0$ (depending on $\|\omega_0\|_{H^2}$ and $\nu$) and  a unique solution $\omega$ of \eqref{eqvorticity3d}-\eqref{BCvorticity3d} such that
\begin{equation} \label{H2_estimate_omega}
     \sup_{0 < t \le T_\nu} \|\omega(t)\|^2_{H^2} + \nu \int_0^T \|\mathbb{P} \omega(s)\|_{H^3}^2 ds \le C_{0,\nu}, \quad \omega\vert_{t=0} = \omega_0.
\end{equation}
\end{thm}    
\begin{rem}
As in dimension $d=2$, existence of weak Leray solutions (starting from $u_0 \in L^2$) is an open problem, {\it cf.} Remark \ref{rem_Leray}. 
\end{rem}
\begin{rem}
Contrary to the two-dimensional case, local well-posedness for data $\omega_0 \in H^1$,is  open. Because the normal component of $\omega$ does not vanish at the boundary, $\omega$ differs from its Leray projection by a harmonic field.
As explained later, this yields a  limitation which comes from the fact that diffusion does not act on $(I-\mathbb{P})\omega$, with $\mathbb{P}$ the Leray projector. This part of $\omega$ somehow obeys an hyperbolic evolution, requiring more regularity for local well-posedness. 
\end{rem}
\begin{rem}
In contrast with the two-dimensional case, the time of existence given in Theorem \ref{thmWP3d} now depends on $\nu$. An interesting open problem is to get existence on a time interval independent of $\nu$, see \cite{MasRou} in the case of Navier boundary condition. 
\end{rem}
\begin{rem} 
The condition \eqref{BCvorticity3d} has some connection with models for perfect conductors in electromagnetism, where  the time harmonic Maxwell equations 
$$ \curl E - i k B = 0, \quad \curl B + ik E = 0 \quad \text{ in } \: \Omega  $$
are completed by 
$$ B \cdot n\vert_{\pa \Omega} = 0, \quad E \times n\vert_{\pa \Omega} = 0. $$
If $B$ is thought as an analogue of $u$ and $E = \frac{i}{k} \curl B$ is thought as an analogue of $\omega$, these last boundary conditions are very close to \eqref{BCvorticity3d}.   
\end{rem}

\subsection{Boundary layer analysis in a half-plane}
Beyond the satisfactory well-posedness theory of \eqref{NS}-\eqref{BC} (or of its equivalent vorticity formulation), the main advantage of the diffusion-free boundary condition \eqref{BC} is to tame the effect of the diffusive term near the boundary. This can be seen rigorously from a boundary layer analysis. We restrict for simplicity to the case of a half-plane $\Omega = \R^2_+$. We show that the boundary layer induced by \eqref{BC2D} has size $\sqrt{\nu}$, but amplitude $\nu$ only. Our main result, to be proved in Section \ref{secBL}, is: 
\begin{thm} \label{thm_BL}
Let $\Omega = \R^2_+$. Let $\omega_0 \in C^\infty_c(\overline{\Omega})$. For any $N > 0$, there exists
\begin{itemize}
\item a velocity field of boundary layer form: 
\begin{equation} \label{Ansatz}
u^{app}(t,x,y) = \sum_{i=0}^N \sqrt{\nu}^i u^i(t,x,y) + \sum_{i=2}^N \sqrt{\nu}^i u^{i,bl}(t,x,y/\sqrt{\nu})     
\end{equation} 
\item a solution $\omega \in L^\infty_{loc}(\R_+ ; L^2(\Omega) \cap L^1(\Omega))$ of \eqref{eqvorticity2d}-\eqref{BCvorticity2d} 
\end{itemize}
such that $\omega\vert_{t=0} = \curl u^{app}\vert_{t=0} = \omega_0$ and : 
$\forall T > 0, \quad \sup_{t\in [0,T]} \|\omega - \curl u^{app}\|_{L^2} \le C_T \nu^{\frac{N}{2}-1} $
\end{thm}
\begin{rem}
The Ansatz \eqref{Ansatz} reflects a boundary layer of size $\sqrt{\nu}$, and as the last sum at the r.h.s. of \eqref{Ansatz} starts at $i=2$, the velocity in the boundary layer has amplitude $\nu$. This is in sharp contrast with the usual no-slip condition (Prandtl layer of amplitude $1$) and stress-free or Lions' boundary conditions (amplitude $\sqrt{\nu}$). 
\end{rem}

\section{2D Well-posedness} \label{secWP2d}
In all what follows, $\| \, \| = \| \, \|_{L^2}$. 
We focus here on the proof of Theorem \ref{thmWP2d}. We follow a classical approach, first we provide {\it a priori} estimates, and then we explain how to use  these formal bounds to build a unique solution. 
\subsection{Estimates for $L^2$ initial vorticity} \label{formal_L2}
Assume that we are given a smooth enough solution $\omega$ of \eqref{eqvorticity2d}-\eqref{BCvorticity2d} with condition 
$\omega\vert_{t=0} = \omega_0$. Multiplying by $\omega$ and integrating by parts, we get 
$$ \frac12 \pa_t \|\omega\|^2 + \nu \| \na \omega\|^2 = - \frac{1}{2} \int_{\pa \Omega} (u \cdot n) \omega^2 - \nu  \int_{\pa \Omega} \pa_n \omega \, \omega = 0 \, , $$
resulting in
\begin{equation} \label{L2_estimate_omega}
\frac12 \|\omega(t)\|^2  + \nu \int_0^t \| \na \omega\|^2 =  \frac12 \|\omega_0\|^2 \, .
\end{equation} 
It is important to stress again that $\Omega$ is bounded and simply connected. We can decompose $\omega = \bar\omega + \tilde\omega$, where $ \bar\omega = \frac{1}{|\Omega|} \int_\Omega \omega$ is the average over the domain. One has 
$$ \pa_t  \bar\omega = - \int_\Omega \div (u \omega)  - \nu \int_\Omega \div (\na \omega) = 0 \, , $$
so that this average  $\bar\omega$ is constant in time
$$ \bar{\omega}(t) = \bar{\omega}_0 \, .$$
For the mean-free part $\tilde\omega$, we can use the Poincaré inequality, and find
\begin{equation} \label{estim_meanfree_vorticity}
\|\tilde \omega(t)\|^2  + \nu \int_0^t \| \tilde\omega\|_{H^1}^2 \le  C \|\omega_0\|^2 \, .
\end{equation}
\begin{rem} \label{rem_exp_decay}
The inequality \eqref{estim_meanfree_vorticity}  yields in particular an exponential decay $\|\tilde\omega(t)\| \le C  e^{- c \nu t}$. 
Moreover,  from the Biot and Savart law $u = \na^\perp \Delta^{-1} \omega$, where $\Delta^{-1}$ refers to the inverse of the  Laplacian (subject to  Dirichlet condition at $\pa \Omega$), we find that 
$$ \|u(t)\|_{H^1} \le C \|\omega\| $$
(more generally, for any integer $m$ such that the r.h.s. is finite,
$ \|u(t)\|_{H^{m+1}}  \le C_k \|\omega(t)\|_{H^m}$).
Following the above decomposition of $\omega$, we can split the stream function $\psi$, satisfying 
$$\Delta \psi = \omega, \quad u = \na^\perp \psi$$ 
as 
$$ \psi = \bar\omega \, \bar\Psi + \tilde \psi $$
where $\Bar\Psi = \Bar\Psi(x,y)$ is the (constant in time) function satisfying 
$$ \Delta  \Bar\Psi  = 1, \quad \bar\Psi_{\mid \pa \Omega} = 0 $$
and where 
$$ \Delta  \tilde \psi = \tilde \omega, \quad \tilde\psi_{\mid \pa \Omega} = 0 \, .$$
We can now decompose $u =  \bar\omega \, \bar U + \tilde u$ with $\bar U = \na^\perp \bar\Psi$, $\tilde u = \na^\perp \tilde\psi$, and 
$$ \|\tilde u(t)\|_{H^1(\Omega)} \le C \|\tilde \omega(t)\| \le C'  e^{- c \nu t}$$
\end{rem}
\noindent To obtain an {\it a priori} estimate for ${d\omega}/{dt}$, we write  
$$ \left\|\frac{d\omega}{dt}\right\|_{(H^1)'} \le \|u \cdot \na \omega\|_{(H^1)'} + \nu \|\Delta \omega\|_{(H^1)'} \, .$$
Because of to the non-penetration condition, we  have for all $\varphi \in H^1$
$$ \langle u \cdot \na \omega , \varphi \rangle_{\langle (H^1)' , H^1 \rangle} =  \langle \div( u \omega) , \varphi \rangle_{\langle (H^1)' , H^1 \rangle} = - \int_\Omega u \omega \cdot \na \varphi \, ,$$
so that for any $\delta > 0$, denoting $\theta = \frac{\delta}{2+\delta}$,
\begin{align*}
\|u \cdot \na \omega\|_{(H^1)'} & \le C \|u \omega\| \le C \|u\|_{L^\frac{2(2+\delta)}{\delta}} \|\omega\|_{L^{2+\delta}} \\
& \le C' \|u\|_{H^1} \|\omega\|^{1 - \theta} \|\omega\|_{H^1}^{\theta} \le C'' \|\omega\|^{2-\theta} \|\omega\|_{H^1}^\theta.     
\end{align*}
Finally, thanks to the Neumann condition on $\omega$, we have for all $\varphi \in H^1$
$$ \langle \Delta \omega , \varphi \rangle_{\langle (H^1)' , H^1 \rangle} = - \int_\Omega \na \omega \cdot \na \varphi $$
so that $\|\Delta \omega\|_{(H^1)'} \le C \|\na \omega\|$. We end up with 
$$ \|\frac{d\omega}{dt}\|_{(H^1)'}  \le C \big( \|\omega\|^{2-\theta} \|\omega\|^{\theta}_{H^1} + \nu \|\omega\|_{H^1} \big) $$
which implies that 
$$\nu^\theta \int_0^T \|\frac{d\omega}{dt}\|_{(H^1)'}^2 \le C \Big( \sup_{0 \le t \le T} \|\omega(t)\|^{4-2\theta} \big(\nu \int_0^T \|\omega\|_{H^1}^2 \big)^\theta +  \nu^{\theta+2} \int_0^T \|\omega\|_{H^1}^2 \Big) \le C'  \, .$$ 
The estimate \eqref{estim2d1} follows from the above observations. 

\subsection{Estimates for $H^1$ initial vorticity.} \label{formalH1}
In case $\omega_0 \in H^1(\Omega)$, we can further derive the estimate \eqref{estim2d2}, as follows: we multiply \eqref{eqvorticity2d} by $\pa_t \omega$ and integrate by parts to find 
\begin{align*}
    \nu \frac12 \frac{d}{dt} \|\na \omega(t)\|^2 +  \|\pa_t \omega\|^2 = - \int_\Omega u \cdot \na \omega \, \pa_t \omega  & \le \frac{1}{4}  \|\pa_t \omega\|^2  +  \|u \cdot \na \omega\|^2 \\
      & \le \frac{1}{4}  \|\pa_t \omega\|^2 +  \|u\|_{L^\infty}^2  \|\na \omega\|^2 \\
      & \le \frac{1}{4}  \|\pa_t \omega\|^2 +  C \|u\|_{H^2}^2  \|\na \omega\|^2  \\
      & \le \frac{1}{4}  \|\pa_t \omega\|^2 +  C' \|\omega\|^2_{H^1}  \|\na \omega\|^2 \, .
    \end{align*}
We thus end up with 
$$   \nu \|\na \omega(t)\|^2 +  \int_0^t \|\pa_t \omega\|^2 \le C \int_0^t \|\omega\|^2_{H^1}  \|\na \omega\|^2 \, ,$$
which by Gronwall's lemma implies 
 $$ \|\na \omega(t)\|^2 \le \|\na \omega(0)\|^2 e^{C/\nu \int_0^t \|\omega\|^2_{H^1}} \le \|\na \omega(0)\|^2  e^{C/\nu^2 \|\omega_0\|^2} \, ,$$
where the last bound stems from \eqref{estim2d1}. Decomposing again $\omega = \bar{\omega} + \tilde{\omega} = \bar{\omega}_0 + \tilde{\omega}$ and using Poincaré inequality for the last part, 
$$ \sup_{0 \le t \le T}\|\omega(t)\|_{H^1}^2 \le C_{0,\nu} \, ,$$
and then 
$$ \int_0^T \|\pa_t \omega\|^2 \le C_{0,\nu} \, . $$ 
To obtain an $H^2$ estimate, we now write 
$$ - \nu \Delta \tilde \omega = - u \cdot \na \omega - \pa_t \omega  \, .$$
Since $\tilde \omega$ has zero mean and satisfies a Neumann condition at $\pa \Omega$, classical elliptic regularity yields 
$$  \|\tilde \omega\|_{H^2} \le \mathcal{C} \big( \|u \cdot \na \omega\| + \|\pa_t \omega\|\big) \, . $$
The second term is controlled by the previous estimate. For the first one, we make use of the  interpolation inequality (valid in dimension 2) : $\|f\|_{L^4} \le C \|f\|^{\frac12} \|f\|_{H^1}^{\frac12}$ to get 
\begin{align*}
    \|u \cdot \na \omega\| & \le \|u\|_{L^4} \|\na \omega\|_{L^4}   \le C  \|u\|_{H^1} \|\na \omega\|^{\frac12} \|\na \omega\|_{H^1}^{\frac12} \\
    & \le C' \|\omega\| \|\omega\|_{H^1}^{\frac12} \|\tilde \omega\|_{H^2}^{\frac12}  \le \frac{1}{4\mathcal{C}} \|\tilde \omega\|_{H^2}  + C''  \|\omega\|^2 \|\omega\|_{H^1} \, .
\end{align*} 
We deduce that $\int_0^t \|\tilde \omega\|_{H^2}^2 \le C_{0,\nu}$. Estimate \eqref{estim2d2} follows from  these bounds. 
\subsection{Stability estimate} \label{formalstability}
We conclude these {\it a priori} estimates with a stability estimate that explains the uniqueness of the solution in Theorem \ref{thmWP2d}. 
Given two smooth enough solutions $\omega_1$ and $\omega_2$ of \eqref{eqvorticity2d}, and introducing $\omega = \omega_1 - \omega_2$, one can check that 
$$ \pa_t \omega + u_1  \cdot \na \omega + u \cdot  \na \omega_2 - \nu \Delta \omega = 0, \quad u = \mathcal{BS} \omega $$
with a Neumann condition on $\omega$. An energy estimate gives 
\begin{align*}
 \frac12 \pa_t \|\omega\|^2 + \nu \| \na \omega\|^2 &  \le \int_{\omega } |u| \, | \na \omega_2| \, |\omega| \le \|\na\omega_2\| \, \|u\|_{L^4} \|\omega\|_{L^4} \\
 & \le C  \|\na \omega_2\|  \, \|\omega\|^{3/2}   \| \omega\|_{H^1}^{1/2}\\
 & \le C' \|\na \omega_2\|  \, \|\omega\|^{2} + C  \|\na \omega_2\|  \, \|\omega\|^{3/2}   \|\na  \omega\|^{1/2} \\
 & \le C'   \|\na \omega_2\|  \, \|\omega\|^{2}  +  \frac{\nu}{2} \| \na \omega\|^2 + C_\nu \|\na \omega_2\|^{4/3}  \, \|\omega\|^{2} \, .
\end{align*}
The last inequality follows from the Young inequality $|ab| \le \frac{\nu}{2} a^{4} + c_\nu b^{4/3}$. By Gronwall inequality, it implies that for some $C > 0$:  
 $$ \|\omega(t)\| \le  \|\omega(0)\| e^{C \int_0^t   (\|\na \omega_2^\nu\| +   \|\na \omega_2^\nu\|^{4/3})}  $$
In particular, if $\omega_1(0) = \omega_2(0)$, we find that $\omega_1 = \omega_2$ for all time, i.e. uniqueness of the solution. 

\subsection{Existence and uniqueness of the solution} \label{rigorous2D}
On the basis of the previous {\it a priori} estimates, one can construct a solution using a Galerkin scheme. First, we need to associate a variational formulation to \eqref{eqvorticity2d}-\eqref{BCvorticity2d}. Given $\omega_0 \in L^2(\Omega)$, and $T > 0$ we want to find 
$$ \omega \in L^\infty(0,T; L^2(\Omega)) \cap L^2(0,T; H^1(\Omega)), \quad \text{ with } \frac{d\omega}{dt} \in L^2(0,T; (H^1(\Omega))')$$
such that: for almost every  $t \in (0,T)$, for all $\varphi \in H^1(\Omega)$ 
\begin{align*}
    \langle \frac{d \omega}{dt} , \varphi \rangle_{\langle H^{-1}, H^1 \rangle}  & + \int_\Omega (u \cdot \na \omega) \varphi  + \nu \int_\Omega \na \omega \cdot \na \varphi = 0, \quad u = \mathcal{BS}\omega, \quad \omega(0) = \omega_0 \, . 
    \end{align*}
Note that the Neumann boundary condition is implicit in the last term at the left-hand side, which corresponds to the formal integration by parts
$$ \int_\Omega (\Delta \omega) \varphi = - \int_\Omega \na \omega \cdot \na \varphi + \int_{\pa \Omega} \pa_n \omega \varphi = - \int_\Omega \na \omega \cdot \na \varphi \, .$$
Note also that the term $\int_\Omega (u \cdot \na \omega) \varphi$ is well-defined. Indeed, by property of the Bio-Savart operator and Sobolev embedding,  $u \in L^2(0,T; L^\infty)$, while $\na \omega \in L^2(0,T; L^2)$ and $\varphi \in L^\infty(0,T; L^2)$.  
Finally, as  $\omega \in L^2(0,T;H^1$ and ${d\omega}/{dt} \in L^2(0,T; (H^1)')$, then $\omega \in C([0,T]; L^2)$: see Theorem II.5.14 in \cite{Boyer}. This gives a meaning to the initial condition $\omega(0) = \omega_0$. 

To solve this variational formulation, we rely on a very classical Galerkin scheme. We introduce an orthonormal basis $(e_n)_{n \in \N}$ of $L^2(\Omega)$, made of eigenfunctions of the Neumann Laplacian. More precisely, we take $e_0 = |\Omega|^{-\frac12}$ and $(e_i)_{i \ge 1}$ an orthonormal basis of 
$$L^2_0(\Omega) = \{u \in L^2(\Omega), \int_\Omega u = 0\}$$
made of eigenfunctions of the Laplacian: they are solving for all $i \ge 1$
$$-\Delta e_i = \lambda_i e_i, \quad \pa_n e_i\vert_{\pa \Omega} = 0, \quad \int_\Omega e_i = 0 $$
for an increasing sequence of positive eigenvalues $(\lambda_i)_{i \ge 1}$. As $\Omega \in C^{1,1}$, $e_i \in H^2(\Omega)$ for all $i \ge 1$. Moreover, $(e_n)_{n \in \N}$ is an orthogonal basis of $H^1(\Omega)$. Let $H_n = \text{span}(e_0, \dots, e_{n-1})$ and $P_n$ the orthogonal projection on $H_n$ (with respect to the $L^2$ scalar product). We look at the sequence of approximate problems: find $\omega_n \in C^1([0,T); H_n)$ satisfying on $(0,T)$, for all $\varphi_n \in H_n$,
\begin{align*}
    \langle \frac{d \omega_n}{dt} , \varphi_n \rangle_{\langle H^{-1}, H^1 \rangle}  & + \int_\Omega (u_n \cdot \na \omega_n) \varphi_n  + \nu \int_\Omega \na \omega_n \cdot \na \varphi_n = 0, \quad u_n = \mathcal{BS}\omega_n, \quad \omega_n(0)  =  P_n \omega_0. 
    \end{align*}
It is well-known that, writing $\omega_n(t) = \sum_{i=0}^{n-1} \omega_{n,i}(t) e_i$, this approximate problem can be reformulated as a system of nonlinear odes on $(\omega_{n,i})_{0 \le i \le n-1}$, with smooth vector field. Hence, one has by the Cauchy-Lipschitz theorem the existence and uniqueness of a solution $\omega_n$ on a maximal time interval $[0,T_n)$, $T_n \le T$. We can then take $\varphi_n = \omega_n$ as a test function, and integrate in time between $0$ and some  $t < T_n$. Relying on Theorem II.5.13 in \cite{Boyer}, we get an estimate like \eqref{L2_estimate_omega}, with $\omega_n$ replacing $\omega$ (note that $\|P_n\omega_0\| \le \|\omega_0\|$). This uniform bound shows that $T_n = T$, and allows to extract a subsequence of $(\omega_n)_n$ that converges weakly* in $L^\infty(0,T; L^2)$ and weakly in $L^2(0,T; H^1)$ to some $\omega$. The energy bound \eqref{estim2d1}, satisfied uniformly in $n$ by $\omega_n$ remains true for $\omega$. Also, for an arbitrary $\varphi \in H^1(\Omega)$, taking $\varphi_n = P_n \varphi$, we get 
\begin{align*}
   \Big|  \langle \frac{d \omega_n}{dt} ,  \varphi  \rangle_{\langle H^{-1}, H^1 \rangle} \Big| & =  \Big|  \int_\Omega \frac{d \omega_n}{dt}   \varphi \Big| =    \Big|  \int_\Omega \frac{d \omega_n}{dt}  P_n \varphi \Big| \\
   & =   \Big| \int_\Omega (u_n \cdot \na \omega_n) P_n \varphi  + \nu \int_\Omega \na \omega_n \cdot \na P_n \varphi  \Big| \\
   & \le \|u_n\|_{L^4} \|\na \omega_n\| \|P_n \varphi\|_{L^4} + \nu \|\na \omega_n\| \|P_n \varphi\|_{H^1}  \\
   & \le C \left(  \|\omega_n\|_{L^2} \|\omega_n\|_{H^1} + \|\omega_n\|_{H^1} \right) \|P_n \varphi\|_{H^1} \\
   & \le C \left(  \|\omega_n\|_{L^2} \|\omega_n\|_{H^1} + \|\omega_n\|_{H^1} \right) \|\varphi\|_{H^1}  \, .
   \end{align*}
   We have used again here the inequalities $\|u_n\|_{H^1} \le C \|\omega_n\|$ (property of the Biot-Savart kernel), $\|P_n \varphi\|_{L^4} \le C \|P_n \varphi\|_{H^1}$ (Sobolev embedding), and the inequality $\|P_n \varphi\|_{H^1} \le \|\varphi\|_{H^1}$. We deduce the pointwise in time bound
   $$\|\frac{d \omega_n}{dt}\|_{(H^1)'} \le C \left(  \|\omega_n\|_{L^2} \|\omega_n\|_{H^1} + \|\omega_n\|_{H^1} \right)  \, ,$$
which shows that ${d \omega_n}/{dt}$ is uniformly bounded in $L^2(0,T; (H^1)')$ and so weakly converges (after extraction) to ${d\omega}/{dt}$ in this space. Standard compactness arguments of Aubin-Lions type (see \cite{Simon}) yield that 
$$ u_n \cdot \na \omega_n  \rightarrow u \cdot \na \omega, \quad u = \mathcal{BS} \omega$$
weakly in $L^2(0,T ; L^{\frac43})$. Taking a fixed $N$, a fixed $\varphi \in H_N$ and letting $n \rightarrow +\infty$, we get for almost every $t$
\begin{align*}
    \langle \frac{d \omega}{dt} , \varphi \rangle_{\langle H^{-1}, H^1 \rangle}  & + \int_\Omega (u \cdot \na \omega) \varphi  + \nu \int_\Omega \na \omega \cdot \na \varphi = 0, \quad u = \mathcal{BS}\omega, \quad \omega(0)  =   \omega_0. 
    \end{align*}
Finally, by density of $\cup_N H_N$ in $H^1$, this equality holds for all $\varphi \in H^1$. This shows that $\omega$ solves the variational formulation. 

Regarding uniqueness of such solution, one can mimic the derivation of the stability estimate above : given two solutions $\omega_1$ and $\omega_2$, one substracts the variational formulation for $\omega_2$ to the one for $\omega_1$, and test with $\varphi = \omega_1 - \omega_2$. The rest of the reasoning is identical. 

\subsection{Extra regularity} \label{rigorous_extra}
We still need to show the last regularity statement of Theorem \ref{thmWP2d}, when $\omega_0 \in H^1(\Omega)$. It is better to go back to the approximate problem for $\omega_n$, and test this approximate variational formulation with $\varphi_n = \pa_t \omega_n$. We can proceed as we did for the {\it a priori} estimate and recover in this way that ${d\omega_n}/{dt}$ is bounded uniformly in $n$ in $L^2(0,T; L^2)$, so that $\|{d\omega}/{dt}\|_{L^2((0,T) \times \Omega)} \le C_{0,\nu}$. To recover the $H^2$ estimate, we again write 
$$ -\nu \Delta \tilde \omega = - u \cdot \na \omega - \pa_t \omega \, .$$
By the previous estimate, for a.e. $t$, $\|\frac{d \omega}{dt}(t)\| < \infty$. Also, 
\begin{align*}
    \|u \cdot \na \omega(t)\| & \le \|u(t)\|_{L^\infty} \|\na \omega(t)\| \le C \|u(t)\|_{H^2} \|\na \omega(t)\| \\
    & \le C' \|\omega(t)\|_{H^1}^2 < \infty \quad \text{ for almost every $t$.} 
\end{align*} 
It follows from elliptic regularity that $\tilde\omega(t)$ belongs to $H^2(\Omega)$ for almost every $t$. From there, one can apply the {\it a priori estimates} seen earlier and conclude that $\int_0^T \|\omega\|_{H^2}^2 \le C_{0,\nu}$. 

\subsection{Extension to the half-plane case} \label{halfplane}
In the case $\Omega = \R^2_+$, the well-posedness statement and its proof have to be slightly modified. The reason is that the Biot-Savart operator $\mathcal{BS} = \na \Delta^{-1}$ associated to the half-plane only obeys the homogeneous estimate  
$$ \|\na \mathcal{BS}f \|_{W^{m,p}} \le C_p \|f\|_{W^{m,p}}, \quad \forall 1 < p < \infty, \quad \forall m \in \N \, ,$$
that is $\|\mathcal{BS}f\|_{W^{m+1,p}}$ is replaced by $\|\na \mathcal{BS}f\|_{W^{m,p}}$. If we stick to the case $p=2$ considered in Paragraphs \ref{formal_L2} and \ref{formalH1}, no direct control of $u$ itself by $\omega$ is available, and some of the estimates fail. To  obtain a control of $u$, we need to consider some $p < 2$, to benefit from Sobolev inequality   
$$ \|\mathcal{BS}f\|_{L^{\frac{2p}{2-p}}} \le C_p \|f\|_{L^p}, \quad \forall 1 \le p < 2.$$
 But of course, this requires an additional control of $\omega(t)$ in $L^p(\R^2_+)$. This control comes from the following {\it a priori} estimate: multiplying \eqref{eqvorticity2d} by $f'(\omega)$ for $f$ a smooth convex function, we get 
$$ \pa_t f(\omega) + u \cdot \na f(\omega) - \nu f'(\omega) \Delta \omega = 0,  $$
and after integration over $\Omega$
\begin{equation} \label{estim_Casimir}
\pa_t \int_\Omega f(\omega) = - \nu \int_\Omega f''(\omega) |\na \omega|^2 \le 0.
\end{equation}
Through a standard approximation argument, we can extend this inequality to $f(x) = |x|^p$ for any $p \ge 1$, resulting in 
$$ \pa_t \int_\Omega |\omega|^p \le 0.$$
We can consider for instance $p=1$: uder the additional  assumption that $\omega_0 \in L^1(\R^2_+)$, we get a control of $\omega$ in $L^\infty(L^1)$ and in turn a control of $u = \mathcal{BS}\omega$ in $L^\infty L^2$. From there, the reasoning of the previous paragraphs can easily be adapted and results in 
\begin{thm} \label{thmhalfplane}   
 Let $\Omega = \R^2_+$. Let $T > 0$. Let $\omega_0 \in L^2(\Omega) \cap L^1(\Omega)$. There exists a unique solution $\omega$ of \eqref{eqvorticity2d}-\eqref{BCvorticity2d} such that, for some $C > 0$ independent of $\nu \le 1$
\begin{equation} 
\sup_{0 < t \le T} \|\omega(t)\|^2_{L^2 \cap L^1} +  \int_0^T \|\frac{d\omega}{dt}(s)\|^2_{(H^1)'} ds +  \nu \int_0^T \|\omega(s)\|_{H^1}^2 ds \le C \|\omega_0\|_{L^2}^2, \quad \omega\vert_{t=0} = \omega_0. 
\end{equation}
Moreover, in the case $\omega_0 \in H^1(\Omega)$, we get the following additional estimate,  for some $C_{0,\nu}$ (depending on $\|\omega_0\|_{H^1}$ and $\nu$)
\begin{equation}
\sup_{0 < t \le T} \|\omega(t)\|^2_{H^1 \cap L^1}  +  \int_0^T \|\frac{d\omega}{dt}(s)\|^2_{L^2} ds + \nu \int_0^T \|\omega(s)\|_{H^2}^2 ds \le C_{0,\nu}.
\end{equation}
\end{thm}    

\subsection{Extension to non-simply connected domains} \label{nonconnected}
Besides the case of a simply connected domain covered by Theorem \ref{thmWP2d}, it is worth mentioning the case of a domain $\Omega$ with holes. Namely, we consider here a 2D domain $\Omega$ of the form 
$$ \Omega = \Omega_0 - \cup_{i=1}^N \overline{\Omega_i}$$ 
where $\Omega_i$ is a smooth bounded simply connected domain for all $0 \le i \le N$ and such that $\overline{\Omega_i} \subset \Omega_0$  for all $1 \le i \le N$. In this setting, it is well-known that special attention must be paid to the construction of the Biot-Savart operator, that is to the reconstruction of $u$ from $\omega$. One can as in the simply connected case try to look for $u$ in the form $u = \na^\perp \psi$, with a stream function $\psi$  satisfying 
$$ \Delta \psi = \omega \: \text{ in } \: \Omega, \quad \pa_\tau \psi = 0  \: \text{ at } \:  \pa\Omega.$$
Note that $\psi$ is defined up to an additive constant, and that the vanishing of its tangential derivative at $\pa \Omega$ corresponds to the non-penetration condition on $u$. In the simply connected case, one can fix the additive constant by enforcing  the Dirichlet condition $\psi\vert_{\pa \Omega} = 0$, which defines a well-posed problem for $\psi$ and allows to recover $u$ from $\omega$.  But in the case of multiple connected components considered here, the additive constant is fixed by the mere condition $\psi\vert_{\pa \Omega^0} = 0$, while the constant values of $\psi$ at the other boundary parts $\pa \Omega_i$, $i \ge 1$, remain {\it a priori} undetermined. Moreover, there is no canonical way to select the values of these constants. Still, one has the following crucial result, that allows to determine these constants, and therefore $\psi$ and $u$, if we know the circulation of $u$ around each obstacle $\Omega_i$, $i \ge 1$. 
\begin{thm}
    For each $1\le i \le N$, there exists a unique (smooth) $\Psi_i$ satisfying 
    $$ \Delta \Psi _i = 0  \: \text{ in } \: \Omega, \quad  \Psi_i\vert_{\pa \Omega_0} = 0, \quad  \pa_\tau \Psi_i\vert_{\pa \Omega_j} = 0, \quad \oint_{\pa \Omega_j} \pa_n \Psi_i = \delta_{ij} \quad \forall j \ge 1.$$
    Moreover,  $\na^\perp \Psi_1, \dots, \na^\perp\Psi_N$ form a basis of the vector fields that are curl-free, divergence-free and tangent at  $\pa \Omega$. It follows that given $\Gamma = (\Gamma_1, \dots, \Gamma_N) \in \R^N$ and a smooth $\omega$,  there exists a unique smooth field $u$ satisfying 
    $$ \curl u = \omega, \quad \div u = 0 \: \text{ in } \: \Omega, \quad u \cdot n\vert_{\pa \Omega} = 0, \quad \oint_{\pa \Omega_i} u \cdot \tau = \Gamma_i, \: i \ge 1.  $$
\end{thm}
\noindent We refer to \cite{Kikuchi,Lopes} for all details. We denote $u = \mathcal{BS}_\Gamma \omega$. It is easily seen that 
$$ u = \na^\perp \psi^0 + \sum_{i=1}^N c_i \na^\perp \Psi_i$$
where $\psi_0$ satisfies 
$$ \Delta \psi^0 = \omega \: \text{ in } \: \Omega, \quad \psi^0\vert_{\pa \Omega} = 0 ,$$
while 
$$ \forall 1 \le i \le N, \quad c_i = \Gamma_i - \oint_{\pa \Omega_i} \na^\perp \psi^0 \cdot \tau = \Gamma_i - \oint \pa_n \psi_0.$$
One can then extend the definition of $\mathcal{BS}_\Gamma \omega$  to vorticities $\omega$ in $L^p$. 

Back to the analysis of \eqref{eqvorticity2d}-\eqref{BCvorticity2d}, we see that we  need to know the circulation of $u$ around each obstacle. In the Euler case, such circulations are preserved through time, by Bernoulli's conservation theorem. Nicely, this is still the case under the diffusion-free boundary condition 
$$ \Delta u \cdot \tau\vert_{\pa \Omega} = 0, \quad  u \cdot n\vert_{\pa \Omega} = 0 .$$
Indeed, this can be seen by coming back to the Navier-Stokes equation in the velocity form, written as  
$$ \pa_t u + \omega u^\perp + \na (\frac{1}{2}|u|^2 + p) - \nu \Delta u = 0.$$
Taking the scalar product with $\tau$ and integrating along each boundary part $\pa \Omega_i$, $i \ge 1$, we get (the diffusion-free boundary condition is crucial here): 
$$ \pa_t \oint_{\pa \Omega_i} u \cdot \tau +  \oint_{\pa \Omega_i} u \cdot n \, \omega + \oint_{\pa \Omega_i} \na (\frac{1}{2}|u|^2 + p) = 0$$
which comes down to $\pa_t \oint_{\pa \Omega_i} u \cdot \tau = 0$. 
Hence, $\Gamma_i(t) = \Gamma_i(0)$ for all $t$, and so $u = \mathcal{BS}_{\Gamma(0)} \omega$. From then on, the mathematical treatment is the same as for the simply connected domain, leading to a similar well-posedness result. 

\section{3D well-posedness} \label{secWP3d}
In this section, we turn to the proof of Theorem \ref{thmWP3d}, analyzing system \eqref{NS}-\eqref{BC}  in its vorticity formulation \eqref{eqvorticity3d}-\eqref{BCvorticity3d}. This last system will reveal significantly harder than its two-dimensional counterpart, due to two factors. First, the stretching term $\omega \cdot \na u$  is a well-known mechanism for growth vorticity, and complicates the energy estimates. Second, contrary to the 2D setting, the condition \eqref{BCvorticity3d} is not equivalent to a Neumann condition on $\omega$. This implies that we do not have the identity  $-\int_\Omega \Delta \omega \cdot \omega = \int_\Omega |\na \omega|^2$, but rather 
\begin{align*}
-\int_\Omega \Delta \omega \cdot \omega  & = \int_\Omega (\curl \curl \omega) \cdot \omega = \int_\Omega |\curl \omega|^2 + \int_{\pa \Omega}   (n \times \curl \omega) \cdot \omega = \int_\Omega |\curl \omega|^2 .
\end{align*}
This yields a control of  $\curl \omega$ in $L^2$, not of the full gradient. It can be tempting at this point to use estimates of Biot-Savart type: following \cite[Remark 1.4 and Lemma 2.15]{Cheng}, for any smooth simply connected domain $\Omega$, for any $s$,  
\begin{equation} \label{BSproperties}
    f \cdot n\vert_{\pa \Omega} = 0 \: \: \text{ or } \: f \times n\vert_{\pa \Omega} = 0 \quad \Rightarrow \quad \|f\|_{H^{s+1}(\Omega)} \le C_s \|\curl f\|_{H^s(\Omega)} 
\end{equation}  
The problem is that we can not take $f=\omega$, as neither its tangential or normal components vanish at the boundary. To be more specific, let us introduce the Leray projector $\mP$ on the space of divergence-free vector fields tangent at $\pa \Omega$. We also define $\mQ = I - \mP$, and write the Liouville decomposition  
\begin{equation} \label{decompo_omega}
\omega = \mP \omega + \mQ \omega = \mP \omega + \na \phi_\omega .
\end{equation}
As $\omega$ and $\mP \omega$ are both divergence-free, and as $\mP \omega \cdot n\vert_{\pa \Omega} = 0$, $\phi_\omega$ satisfies the Neumann problem
\begin{equation} \label{eq_phi}
    \Delta \phi_\omega = 0 \quad \text{in } \Omega, \quad \pa_n \phi_\omega = \omega \cdot n \quad \text{at } \pa \Omega .
\end{equation}
The difficulty is that, while we have the inequalities
$$\|\mP\omega\|_{H^{s+1}(\Omega)}  \le C_s \|\curl \mP\omega\|_{H^s(\Omega)} =\|\curl \omega\|_{H^s(\Omega)} ,$$
no control of $\mQ \omega$ is available. The diffusion simply does not act on $\mQ \omega$: $\Delta \mQ \omega = \curl \curl \na \phi_\omega = 0$. Therefore, to construct solutions, we will need to adopt another strategy than the one in Section \ref{secWP2d}, closer to the ones used for hyperbolic equations, for which local in time solutions can be obtained in $H^s$, $s >\frac{3}{2}$. This explains the higher regularity $\omega_0 \in H^2(\Omega)$ in Theorem \ref{thmWP3d}. However, differentiating the equation is uneasy due to the boundary: we shall rely on the use of tangential derivatives, see Paragraph \ref{subsec_tangential} for a reminder. The derivation of the {\it a priori} estimate in \eqref{H2_estimate_omega} will then be performed in Paragraph \ref{subsec_apriori}. Finally, a description of a construction scheme for the solution of \eqref{eqvorticity3d}-\eqref{BCvorticity3d} will be given in Paragraph \ref{subsec_construction}, concluding the proof of Theorem \ref{thmWP3d}. To find the right approximation scheme, which allows to get solutions by compactness, is not straightforward. Indeed, the {\it a priori} estimates of Subsection \ref{subsec_apriori} use significantly the structure of the original equation, and this structure is not preserved by arbitrary approximations. 
\begin{rem}{(why not a Neumann condition ?)}
It is tempting to get inspiration from our 2D analysis, and  to replace \eqref{BCvorticity3d} by the set of boundary  conditions
$$u \cdot n \vert_{\pa \Omega} = 0, \quad \pa_n \omega\vert_{\pa \Omega} = 0$$
in order to restore the control on $\na \omega$. And indeed, given a divergence-free vector field $u$, one can show that the equation
\begin{equation} \label{omega3d}
    \pa_t \omega + u \cdot \na \omega - \nu \Delta \omega = \omega \cdot \na u
\end{equation}
 is well-posed when endowed with a Neumann condition, starting  from $\omega_0 \in L^2(\Omega)$. One may then think of solving \eqref{eqvorticity3d} iteratively: 
 $$u_n \rightarrow \omega_{n+1} \rightarrow u_{n+1} = \mathcal{BS} \, \omega_{n+1}$$
 But there is an issue: even if the initial datum $\omega_0$ is divergence-free, the corresponding solution $\omega$  does not satisfy $\div \omega = 0$ at positive times, which prevents to apply the Biot-Savart operator. More precisely, one can check that for a given divergence-free $u$ and for $\omega$ solution of \eqref{omega3d}, the divergence $d_\omega = \div \omega$ satisfies 
$$ \pa_t d_\omega + u \cdot \na d_\omega  - \nu \Delta d_\omega = 0   $$
As there is no clear boundary condition on $d_\omega$ at $\pa \Omega$, one can not deduce that $d_\omega$ is zero at positive times even if it is initially. In the case of \eqref{BCvorticity3d}, we will not be facing this problem. We will use the following alternative formulation of \eqref{eqvorticity3d}:
\begin{equation} \label{eqvorticity3dbis}
\pa_t \omega + u \cdot \na \omega + \nu \curl \curl \omega = \omega \cdot \na u, \quad u = \mathbb{BS}\omega 
\end{equation}
using the identity $-\Delta \omega = \curl \curl \omega$ valid for divergence-free vector fields. On one hand, the expression  $\curl \curl \omega$ is adapted to our boundary condition, as seen above. On the other hand,  for a given divergence-free $u$ and $\omega$ solution of \eqref{eqvorticity3dbis}, we find this time 
$$ \pa_t d_\omega + u \cdot \na d_\omega  = 0 $$
so the divergence-free constraint is preserved for all times. See Subsection \ref{subsec_construction} for more on the construction scheme. 
\end{rem}

\subsection{Reminder on tangential derivatives} \label{subsec_tangential}
We remind here of one possible construction  of a family $T_1$, \dots, $T_N$ of tangential derivatives, which will be used to derive  {\it a priori} estimates. The first step is to cover the compact set $\overline{\Omega}$ by a finite number of smooth domains $U_0, \dots, U_m$ satisfying the following properties: $U_0 \Subset \Omega$ and for all $i= 1...m$, there exists a diffeomorphism $\Phi_i = \Phi_i(X,Y,Z)$ from a neighborhood of the unit ball $ B(0,1)$ of $\R^3$ to a neighborhood of $U_i$ satisfying 
$$ \Phi_i(B(0,1))  = U_i,  \quad \Phi_i(B(0,1) \cap \{Z > 0\})= U_i \cap \Omega, \quad \Phi_i(B(0,1) \cap \{Z=0\}) = U_i \cap \pa \Omega$$
Let $(\chi_i)_{0 \le i \le m}$ a smooth partition of unity associated to $(U_i)_{0 \le i \le m}$. We then take 
$$T_1  = \chi_0 \, \pa_x, \quad T_2 = \chi_0 \, \pa_y, \quad T_3 = \chi_0 \, \pa_z $$
and for all $1 \le i \le m$
$$ T_{2+2i}  f= \chi_i \, \pa_X (f \circ \phi_i) \circ \phi_i^{-1}, \quad T_{3+2i} f = \chi_i \pa_Y (f \circ \phi_i) \circ \phi_i^{-1} .$$
We shall also note, for any $N$-tuple $\alpha = (\alpha_1, \dots, \alpha_N)$
$$ T^\alpha f = T^{\alpha_1}_1 \dots  T^{\alpha_N}_N f . $$

\subsection{A priori estimates} \label{subsec_apriori}
We are now ready to perform {\it a priori} estimates on (smooth enough) solutions of \eqref{eqvorticity3d}-\eqref{BCvorticity3d}. We first collect several estimates deduced from \eqref{BSproperties}: first, as $\curl \mP \omega = \curl \omega$ and $\mP \omega \cdot n\vert_{\pa \Omega} = 0$, we get for all $s$:
\begin{equation} \label{estimBS1}
\|\mP \omega\|_{H^{s+1}} \le C \|\curl \omega\|_{H^s}
\end{equation}
Then, as $\Delta \omega = -\curl \curl \omega$ and $\curl\omega \times n\vert_{\pa \Omega} = 0$, 
\begin{equation} \label{estimBS2}
\|\curl \omega\|_{H^{s+1}} \le C \|\Delta \omega\|_{H^s}
\end{equation}
which combined with the previous inequality yields 
\begin{equation} \label{estimBS2bis}
\|\mP \omega\|_{H^{s+2}} \le C \|\Delta \omega\|_{H^s} \, .
\end{equation}
Finally, one can check using curvilinear coordinates that $f \times n\vert_{\pa \Omega} = 0 \: \Rightarrow \: \curl f \cdot n\vert_{\pa \Omega} = 0$. Applying this property with $f=\curl \omega$, we get $\Delta \omega\cdot n\vert_{\pa \Omega} = 0$ and it follows again from  \eqref{BSproperties} that 
\begin{equation} \label{estimBS3}
\|\Delta \omega\|_{H^{s+1}} \le C \|\curl \Delta \omega\|_{H^s}
\end{equation}
which together with \eqref{estimBS2bis} yields
\begin{equation} \label{estimBS3bis}
\|\mP \omega\|_{H^{s+3}} \le C \|\curl \Delta \omega\|_{H^s} = C \|\Delta \curl \omega\|_{H^s} \, .
\end{equation}
We now start the estimates. 
We take the curl of  \eqref{eqvorticity3d} and test against $\curl \Delta  \omega$ to get 
\begin{align*}
    \frac{1}{2} \pa_t \|\Delta \omega\|^2 + \frac{\nu}{2} \|\Delta \curl \omega\|^2 & \le C_\nu \| \curl(u \cdot \na \omega - \omega \cdot \na u) \| \\
    & \le C'_\nu \left( \|u\|_{L^\infty} \|\omega\|_{H^2} + \|\na u\|_{L^4} \|\na \omega\|_{L^4} +  \|\na^2 u\|_{L^4} \|\omega\|_{L^4} \right)^2 \\
    & \le C''_\nu \|\omega\|_{H^2}^4 \, .
\end{align*}
Combining this inequality with \eqref{estimBS2bis}-\eqref{estimBS3bis}, we have
\begin{equation} \label{estimPomega}
    \|\mP\omega(t)\|^2_{H^2} + \int_0^t  \|\mP\omega\|^2_{H^3} \le  C_\nu \|\mP\omega(0)\|^2_{H^2} + C_\nu \int_0^t  \|\omega\|_{H^2}^4.
\end{equation}
We now turn to the estimates of tangential derivatives $T^\alpha \omega$ (see Subsection \ref{subsec_tangential}), for $|\alpha| \le 2$. We have
\begin{align*}
    \pa_t T^\alpha \omega + u \cdot \na T^\alpha \omega + \nu T^\alpha \curl \curl \omega = T^\alpha (\omega \cdot \na u) + [u \cdot \na , T^\alpha] \omega \, .
\end{align*}
We deduce 
\begin{align*}
    \frac{1}{2} \|T^\alpha \omega\|^2 & + \nu \int_\Omega T^\alpha  \curl \curl \omega \cdot T^\alpha \omega  = \int_\Omega T^\alpha (\omega \cdot \na u) \cdot T^\alpha \omega + \int_\Omega [u \cdot \na , T^\alpha] \omega \cdot T^\alpha \omega \, .
\end{align*}
We have 
\begin{align*}
    \Big| \int_\Omega T^\alpha (\omega \cdot \na u) \cdot T^\alpha \omega \Big| & \le C  \Big( \|T^\alpha \omega\| \|\na u\|_{L^\infty} + \|\na \omega\|_{L^4} \|\na^2 u\|_{L^4} + \| T^\alpha\na u\| \|\omega\|_{L^\infty} \Big) \|T^\alpha \omega\| \\
    & \le C' \|\omega\|_{H^2}^2 \|T^\alpha \omega\| \le C'' \|\omega\|_{H^2}^3  
\end{align*}
Similarly 
\begin{align*}
    \Big|  \int_\Omega [u \cdot \na , T^\alpha] \omega \cdot T^\alpha \omega \Big| & \le C  \Big( \|\na^2 u\|_{L^4} \|\na \omega\|_{L^4} +  (\|\na u\|_{L^\infty} + \|u\|_{L^\infty}) \|\omega\|_{H^2} \Big) \|T^\alpha \omega\| \\
    & \le C' \|\omega\|_{H^2}^2 \|T^\alpha \omega\| \le C'' \|\omega\|_{H^2}^3  \, .
\end{align*}
We still need to treat the diffusion term
\begin{align*}
   & \int_\Omega T^\alpha  \curl \curl \omega \cdot T^\alpha \omega  =  \int_\Omega [T^\alpha , \curl ] \curl \omega \cdot T^\alpha \omega + \int_\Omega  \curl T^\alpha   \curl \omega \cdot T^\alpha \omega \\
 = & \int_\Omega [T^\alpha , \curl ] \curl \omega \cdot T^\alpha \omega + \int_{\pa \Omega}  n \times (T^\alpha   \curl \omega) \cdot T^\alpha \omega + \int_\Omega   T^\alpha   \curl \omega \cdot \curl T^\alpha \omega \\
 = & \int_\Omega [T^\alpha , \curl ] \curl \omega \cdot T^\alpha \omega + \int_{\pa \Omega}  n \times (T^\alpha   \curl \omega) \cdot (n \times T^\alpha \omega) + \int_\Omega   [T^\alpha ,  \curl] \omega \cdot \curl T^\alpha \omega + \| \curl T^\alpha \omega \|^2 \\
 = & \: I_1 + I_2 + I_3 + \| \curl T^\alpha \omega \|^2 \, .
    \end{align*}
We have 
\begin{align*}
    |I_1| \le C \|\curl \omega\|_{H^2}\|T^\alpha \omega\| & = C  \|\curl \mP \omega\|_{H^2}\|T^\alpha \omega\| \le C  \|\mP \omega\|_{H^3}\|\omega\|_{H^2} \\
    & \le \frac{C}{2} \|\mP \omega\|_{H^3}^2 + \frac{C}{2} \|\omega\|_{H^2}^2\, .
\end{align*}
Similarly, for arbitrary $\delta > 0$
\begin{align*}
    |I_3| \le C \|\curl T^\alpha \omega\| \|\omega\|_{H^2} \le  
    & \le \delta \|\curl T^\alpha \omega\|^2 + C_\delta \|\omega\|_{H^2}^2 \, .
\end{align*}
Regarding the boundary term, we can notice that $\curl \omega \times n\vert_{\pa \Omega} = 0$, so $T^\alpha \big(\curl \omega \times n \big)\vert_{\pa \Omega} = 0$ and so 
$n \times (T^\alpha   \curl \omega)$ is at $\pa \omega$ a combination of components of $\curl \omega$ and $\na \curl \omega$. It follows that 
\begin{align*}
      |I_2| & \le  \|T^\alpha \curl \omega \times n\|_{H^{1/2}(\pa \Omega)} \, \|T^\alpha \omega \times n\|_{H^{-1/2}(\pa \Omega)} 
      \\ 
      & \le C \|(1,\na) \curl  \omega\|_{H^{1/2}(\pa \Omega)} \, \|T^\alpha \omega \times n\|_{H^{-1/2}(\pa \Omega)} \\
      & \le C' \|\curl \omega\|_{H^2} \big(\|T^\alpha \omega\| + \|\curl T^\alpha \omega\|\big) \\
      & \le C_\delta \|\mP \omega\|_{H^3}^2 + \delta \|\curl T^\alpha \omega\|^2 \, .
 \end{align*}
If we gather all bounds, we end up for $\delta$ small enough with 
\begin{equation} \label{estimTalpha_omega}
\begin{aligned}    
\|T^\alpha \omega(t)\|^2 & + \int_0^t  \|\curl T^\alpha \omega\|^2 \\
&\le C_\nu \|T^\alpha \omega(0)\|^2 + C_\nu \big( \int_0^t \|\mP \omega\|^2_{H^3} +  \int_0^t \|\omega\|^2_{H^2} +  \int_0^t \|\omega\|^3_{H^2} \big) \, .
\end{aligned}
\end{equation}
Combining \eqref{estimPomega} (multiplied by a large constant) with \eqref{estimTalpha_omega} we deduce
\begin{equation} \label{estimalmost_omega}
\begin{aligned}    
& \sum_{|\alpha|\le 2} \|T^\alpha \omega(t)\|^2 + \|\mP\omega(t)\|_{H^2}^2  + \int_0^t  \|\mP\omega\|_{H^3}^2  \\
&\le C_\nu \|\omega(0)\|^2_{H^2} + C_\nu \big( \int_0^t \|\omega\|^2_{H^2} +  \int_0^t \|\omega\|^4_{H^2} \big) \, .
\end{aligned}
\end{equation}
The left-hand side of this inequality almost controls the $H^2$ norm of $\omega$, but some (normal) derivatives of $\mQ \omega$ are missing. We recall that $\mQ \omega  = \na \phi_\omega$ for a potential $\phi_\omega$ solving \eqref{eq_phi}.  To be more explicit, let $X,Y,Z$ local curvilinear coordinates  near an arbitrary point of $\pa \Omega$, with $\pa_X,\pa_Y$ tangent at $\pa \Omega$. We need to control $\pa_X^\alpha \pa_Y^\beta \pa_Z^\gamma \phi_\omega$, where $1 \le \alpha + \beta + \gamma \le 3$ locally in $L^2$, whereas we already control $\pa^{\alpha'}_X \pa^{\beta'}_Y (\pa_X,\pa_Y,\pa_Z) \phi_\omega$, $0 \le \alpha' + \beta' \le 2$. Hence, it remains to control $\pa^2_Z(1, \pa_X,\pa_Y,\pa_Z) \phi_\omega$. Denoting $g$ the euclidean metric in this system of coordinates, we find  
$$ 0 =  \Delta \phi_\omega = g_{ZZ}\pa_Z^2 \phi_\omega + \Delta' \phi_\omega $$
with $\Delta'$ a second order operator involving at most one $Z$ derivative. In particular, 
$(\pa_X,\pa_Y,1) \Delta' \phi$ is locally in $L^2$ and by the equation so is $\pa^2_Z(1, \pa_X,\pa_Y) \phi_\omega$. Differentiating once more the equation in $Z$, we get that $\pa^3_Z \phi$ is locally in $L^2$, resulting in 
$$ \|\mQ \omega\|_{H^2} \le  C \sum_{\alpha \le 2} \|T^\alpha \mQ \omega\| \le C' \big(\sum_{\alpha \le 2} \|T^\alpha \omega\| + \|\mP \omega \|_{H^2}\big)$$
Injecting this into \eqref{estimalmost_omega}, we deduce 
\begin{equation} \label{estimfinal_omega}
\begin{aligned}    
& \|\omega(t)\|_{H^2}^2    + \int_0^t  \|\mP\omega\|_{H^3}^2  \le C_\nu \|\omega(0)\|^2_{H^2} + C_\nu \big( \int_0^t \|\omega\|^2_{H^2} +  \int_0^t \|\omega\|^4_{H^2} \big) \, .
\end{aligned}
\end{equation}
It follows from a standard argument that over some small time $T_\nu$, the bound  \eqref{H2_estimate_omega} holds. 

\subsection{Constructing the solution of \eqref{eqvorticity3d}-\eqref{BCvorticity3d} } \label{subsec_construction}
We now explain how to construct rigorously solutions of \eqref{eqvorticity3d}-\eqref{BCvorticity3d}. The key step will be to solve the linearized problem 
\begin{equation} \label{linear}
    \begin{aligned}
        \pa_t \omega + v \cdot \na \omega - \omega \cdot \na v - \Delta \omega & = 0, \quad \text{ in } \: \Omega \\
        \div \omega & = 0, \quad \text{ in } \: \Omega \\
        \curl \omega \times n\vert_{\pa \Omega} & = 0, \quad \omega\vert_{t=0} = \omega_0  
    \end{aligned}
\end{equation}
for a given divergence-free vector field $v$. Specifically, the point is to show that {\em for all $T > 0$, if $v \in L^\infty(0,T; H^3)$, there exists a solution $\omega \in L^\infty(0,T; H^2)$ s.t. $\mP \omega \in L^2(0,T; H^3)$.}

Let us assume temporarily that we have achieved this task. Then, through a slight modification of the estimates in Subsection \ref{subsec_apriori}, one can deduce that: for all $0 < T < 1$, for all $M > 0$, if $\|v\|_{L^\infty H^3} \le M$ then
\begin{equation*}
    \forall t \in (0,T), \quad \|\omega(t)\|_{H^2}^2 + \int_0^t \|\mP\omega\|^2_{H^3} \le C_0 \|\omega_0\|_{H^2}^2 + C_M \int_0^t \|\omega\|_{H^2}^2 
\end{equation*}
for some  constant $C_0$ only depending on $\nu$ and some constant $C_M$ depending on $M$ and $\nu$. It follows that 
$$ \|\omega(t)\|^2_{H^2} \le  C_0 \|\omega_0\|_{H^2}^2 \, e^{C_M t}  $$
and so 
$$ \|\mathcal{BS}\omega(t)\|_{H^3}^2 \le \|\mathcal{BS}\|_{L(H^2,H^3)}  \sqrt{C_0}  \|\omega_0\|_{H^2} \, e^{\frac{C_M}{2} t}  $$
and so for $M = 2 \|\mathcal{BS}\|_{L(H^2,H^3)}  \sqrt{C_0}  \|\omega_0\|_{H^2}$ and $T$ such that $e^{\frac{C_M}{2} T} < 2$, we find that 
$\|\mathcal{BS}\omega\|_{L^\infty H^3} \le M$. From there, standard arguments allow to solve \eqref{eqvorticity3d}-\eqref{BCvorticity3d} by iteration: 
\begin{equation} \label{iteration}
    \begin{aligned}
        \pa_t \omega^{n+1} + u^n \cdot \na \omega^{n+1} - \omega^{n+1} \cdot \na v - \Delta  \omega^{n+1} & = 0, \quad \text{ in } \: \Omega \\
        \div  \omega^{n+1} & = 0, \quad \text{ in } \: \Omega \\
        \curl  \omega^{n+1} \times n\vert_{\pa \Omega} & = 0, \quad  \omega^{n+1}\vert_{t=0} = \omega_0.  
    \end{aligned}
\end{equation}
and updating the velocity field through $u^{n+1} = \mathcal{BS}\omega^{n+1}$. Moreover, to solve \eqref{iteration}, one can  always consider the limit as $k \rightarrow +\infty$ of smooth approximations $u^{n,k}$ of $u^n$ satisfying 
$$ \|u^{n,k}\|_{L^\infty(H^3)} \le M, \quad u^{n,k} \rightarrow u^n \text{ in } L^2((0,T)\times \Omega) \: \text{ weakly}$$
and similarly smooth approximations $\omega_0^k$ of $\omega_0$, and again invoke standard compactness arguments. This means that we can always assume that $v$ and $\omega_0$ are smooth in system \eqref{linear}.     

\medskip
We now turn to the core part, that is the study of \eqref{linear} for smooth $v$ and $\omega_0$. This is not a so classical system, due to the degeneracy of the diffusion term: as explained earlier, the boundary condition only allows for a control of $\curl \omega$ through diffusion, not of the full gradient. To solve it, we will consider the following approximate systems, parametrized by $0 < \eps \ll 1$: 
\begin{equation} \label{NS_stress_free}
    \begin{aligned}
        \pa_t \omega + v \cdot \na \omega - \omega \cdot \na v + \na p - \Delta \omega - \eps \Delta \omega & = 0, \quad\text{ in } \Omega, \\
          \div \omega & = 0, \quad \text{ in } \: \Omega \\
           \curl \omega \times n\vert_{\pa \Omega} + (\eps \pa_n \omega - p n)\vert_{\pa \Omega} & = 0, \quad \omega\vert_{t=0} = \omega_0  
    \end{aligned}
\end{equation}
The main points of these approximations are that 
\begin{itemize}
    \item the evolution equation in \eqref{linear} is turned into a Navier-Stokes type equation on $\omega$, with the introduction of an extra pressure gradient $\na p$ (and an extra diffusion term $-\eps \Delta \omega$). This ensures that the divergence-free condition on $\omega$ is preserved through time. This is a very important property for the derivation of good estimates, and is the main reason for this choice of approximation. 
    \item a term $\eps \pa_n \omega - q n$ of Neumann type is added in the boundary condition. This kind of Neumann condition, together with Navier-Stokes type equations, is well-known in the context of free-surface viscous problems,  where continuity of the full stress tensor is involved: see for instance \cite[chapter IV.7]{Boyer} for the steady Stokes equation, or \cite{Beale,GuoTice} and references therein for free-surface problems. Let us stress that these last papers treat much more difficult systems than the simple linear model \eqref{NS_stress_free}. 
\end{itemize}
The existence and uniqueness of a smooth solution $\omega$ of \eqref{NS_stress_free}  on an arbitrary time interval $(0,T)$  for fixed $\eps$ can be handled by very classical arguments, as the diffusion term yields then a control of $\|\curl \omega\|^2 + \eps \|\na \omega\|^2$, that is of the full gradient. Of course, the bounds are not {\it a priori} uniform in $\eps$, and the keypoint is now to obtain uniform $L^\infty H^2$ estimates for $\omega^\eps$, resp. uniform $L^\infty H^3$ estimates for $\mP\omega^\eps$. Unfortunately, one can not adapt directly all estimates from Subsection \ref{subsec_apriori}, notably those on $\|\mP \omega\|_{L^2 H^3}$: we could not find a quick alternative to testing the equation against $\curl \Delta \omega$. We proceed in several steps, that we describe briefly. 

\medskip
{\em Step 1.  $L^2 H^2$ estimates for $\mP\omega$.}  The standard estimate yields 
$$ \frac{1}{2}\|\omega(t)\|^2 + \int_0^t \|\curl\omega(t)\|^2 +  \eps \int_0^t \|\na\omega(t)\|^2 \le \frac{1}{2}\|\omega_0\|^2 
$$ 
We then multiply the evolution equation in \eqref{NS_stress_free} with $\pa_t \omega$ and integrate in space time. Standard manipulations yield 
\begin{align}
\int_0^t \|\pa_t \omega\|^2 + \|\curl\omega(t)\|^2 +  \eps \|\na\omega(t)\|^2 & \le 
     C  \Big( \|\omega_0\|_{H^1}^2 +  \int_0^t \|v \cdot \na \omega - \omega \cdot \na v\|^2 \Big) \\
     & \le C \|\omega_0\|_{H^1}^2 + C_v \int_0^t  \|\omega\|_{H^1}^2 \, .
\end{align}
Combining the above inequalities  with the properties of the Biot-Savart operator, we get
\begin{align}
\int_0^t \|\pa_t \omega\|^2 + \|\mP\omega(t)\|_{H^1}^2 +  \eps \|\omega(t)\|_{H^1}^2 &  \le C \|\omega_0\|_{H^1}^2 + C_v \int_0^t  \|\omega\|_{H^1}^2 \, .
\end{align}
We then make use of the following lemma, to be proved in Appendix \ref{appendixA}:
\begin{lemme} \label{lemme_Stokes}
    Let $w,f,g,f_b$ smooth vector fields satisfying 
\begin{equation} \label{stokes_eq}
\begin{aligned}
    -\Delta w - \eps \Delta w + \na q  & = f, \\
    \div \omega & = g, \\
    \curl w \times n\vert_{\pa \Omega} + (\eps \pa_n w - q n)\vert_{\pa \Omega} & = f_b.
\end{aligned}
\end{equation}
Then, for all $k \in \N$
\begin{equation} \label{stokes_regularity_estimate}
    \|\mP w\|_{H^{k+2}} + \eps \|w\|_{H^{k+2}} + \|q\|_{H^{k+1}/\R} \le C_k \left(\|f\|_{H^k} + \|g\|_{H^{k+1}} +  \|f_b\|_{k+\frac12} + \|w\|_{H^{k+1}} \right)
\end{equation}
\end{lemme}
We apply it with $k=0$, $w = \omega$, $q=p$, $f = -\pa_t \omega - v \cdot \na \omega + \omega \cdot \na v$. This results in 
\begin{align*}
\int_0^t \big(  \|\mP\omega\|^2_{H^2} + \eps^2 \|\omega\|_{H^2}^2 + \|p\|^2_{H^1/\R}  \big) & \le C \int_0^t \|\pa_t \omega\|^2 + C  \int_0^t \|v \cdot \na \omega - \omega \cdot \na v\|^2  \\
&  \le C \|\omega_0\|_{H^1}^2 + C  \int_0^t \|\pa_t \omega\|^2  + C_v \int_0^t  \|\omega\|_{H^1}^2 \, .
\end{align*}
Together with the previous estimate, we find 
\begin{align*}
\int_0^t \big( \|\pa_t \omega\|^2 + \|\mP\omega\|^2_{H^2} & + \eps^2 \|\omega\|_{H^2}^2 + \|p\|^2_{H^1/\R} \big)  + \|\mP\omega(t)\|_{H^1}^2 +  \eps \|\omega(t)\|_{H^1}^2 \\
&  \le C \|\omega_0\|_{H^1}^2 + C_v \int_0^t  \|\omega\|_{H^1}^2 \, .
\end{align*}

\medskip
{\em Step 2.  $L^\infty H^1$ estimate for $\omega$.} 
To control the $L^\infty H^1$ norm of $\omega$, we first consider $T\omega$, where $T$ is any tangential derivative, see Subsection \ref{subsec_tangential}. Performing the same type of estimate as in Subsection \ref{subsec_apriori}, we obtain
\begin{align*}
\|T\omega(t)\|^2 + \int_0^t \big(\|T\curl \omega\|^2 & + \eps \|T\na \omega\|^2 \big) \\
& \le C \|T\omega_0\|^2 +   C_v \int_0^t  \|\omega\|_{H^1}^2 + C \int_0^t \big(\|\mP\omega\|^2_{H^2} +  \eps^2 \|\omega\|_{H^2}^2 + \|p\|_{H^1/\R}^2 \big) 
\end{align*}
Finally, we recover the control of the normal derivative of $\omega$ exactly as in Subsection \ref{subsec_apriori}, using that $\mQ\omega = \na \phi_\omega$ with $\phi_\omega$ harmonic. Putting all  estimates together, we end up with 
\begin{equation} \label{uniformH1}
    \|\omega(t)\|^2_{H^1} + \int_0^t \big( \|\mP \omega\|^2_{H^2} + \eps^2 \|\omega\|_{H^2}^2 +  \|p\|_{H^1/\R}^2 \big) \le  C \|\omega_0\|_{H^1}^2 +   C_v \int_0^t  \|\omega\|_{H^1}^2
\end{equation}
which provides an $\eps$-uniform $L^\infty H^1$ bound on $\omega$, resp. $L^2 H^2$ bound on $\mP \omega$. 

\medskip
{\em Step 3. $L^\infty H^2$ bound for $\omega$, $L^2 H^3$ bound for $\mP \omega$} \\
To upgrade the uniform estimate \eqref{uniformH1}, the easiest way (although non optimal) is to differentiate in time  \eqref{NS_stress_free}. The equation satisfied by $\dot{\omega} := \pa_t \omega$ is 
\begin{align*} 
\pa_t \dot{\omega} + v \cdot \na \dot{\omega} - \dot{\omega} \cdot \na v + \na \dot{p} - \Delta \dot{\omega} - \eps \Delta \dot{\omega} &  = - \dot{v} \cdot \na \omega + \omega \cdot \na \dot{v}, \\
\div \dot{\omega} & = 0, \\
 \curl \dot{w} \times n\vert_{\pa \Omega} + (\eps \pa_n \dot{w} - \dot{p} n)\vert_{\pa \Omega} & = 0. 
\end{align*}
We can  proceed with $\dot{\omega}$ as for $\omega$, taking into account the extra source term $- \dot{v} \cdot \na \omega + \omega \cdot \na \dot{v}$: we get the $L^2$ estimate 
\begin{align*} 
\|\dot{\omega}(t)\|_{L^2}^2 & \le  C \|\dot{\omega}\vert_{t=0}\|_{L^2}^2 + C_v \int_0^t \|\omega\|_{H^1}^2 \\
& \le  C \|\omega_0\|_{H^2}^2 + C_v \int_0^t \|\omega\|_{H^1}^2 
\end{align*}
and  eventually
\begin{align*}
    \|\dot{\omega}(t)\|^2_{H^1} + \int_0^t \big(\|\mP \dot{\omega}\|^2_{H^2} + \eps^2 \|\dot{\omega}\|_{H^2}^2 +  \|\dot{p}\|_{H^1/\R}^2 \big) & \le  C \|\pa_t\omega\vert_{t=0}\|_{H^1}^2 +   C_v \int_0^t  \|\dot{\omega}\|_{H^1}^2 + C_v \int_0^t  \|\omega\|_{H^2}^2 \\
    & \le  C \|\omega_0\|_{H^3}^2 +   C_v \int_0^t  \|\dot{\omega}\|_{H^1}^2 + C_v \int_0^t  \|\omega\|_{H^2}^2 
 \end{align*}
where the last integral at the right-hand side is related to the source term. Gronwall's lemma yields 
$$  \|\dot{\omega}(t)\|^2_{H^1} \le C_v \|\omega_0\|_{H^3}^2 + C_v \int_0^t  \|\omega\|_{H^2}^2  \, .  $$
From there, one can go back to Lemma \ref{lemme_Stokes} with  $k=0$ (resp. $k=1$), $w = \omega$, $q=p$, $f = -\pa_t \omega - v \cdot \na \omega + \omega \cdot \na v$, to get an $L^\infty H^2$ bound on $\mP \omega$ (resp. $L^2 H^3$). Precisely,
\begin{align*}
   \|\mP\omega(t)\|^2_{H^2}  + \eps^2 \|\omega(t)\|_{H^2}^2 + \|p\|^2_{H^1/\R}  
& \le C \|\omega_0\|_{H^2}^2 +  C_v \int_0^t \|\omega\|_{H^1}^2 +  C_v   \|\omega(t)\|_{H^1}^2, \\
 \int_0^t  \big( \|\mP\omega\|^2_{H^3}  + \eps^2 \|\omega\|_{H^3}^2 + \|p\|^2_{H^2/\R} \big)  
& \le C_v t \|\omega_0\|_{H^3}^2 + C_v  \int_0^t  \|\omega\|_{H^2}^2 \, .
\end{align*}
Finally, to get the $L^\infty H^2$ bound, we proceed as in Subsection \ref{subsec_apriori}: we perform estimates on $T^\alpha \omega$, $|\alpha| \le 2$, and recover the normal derivatives afterwards. We end up with 
\begin{equation} 
    \|\omega(t)\|^2_{H^2} + \int_0^t \big(  \|\mP \omega\|^2_{H^3} + \eps^2 \|\omega\|_{H^3}^2 +  \|p\|_{H^2/\R}^2 \big) \le  C \|\omega_0\|_{H^3}^2 +   C_v \int_0^t  \|\omega\|_{H^2}^2 + C_v \|\omega(t)\|_{H^1}
\end{equation}
which provides the desired $\eps$-uniform estimate.

\medskip
Thanks to these uniform bounds,  usual compactness arguments provide a sequence of solutions of \eqref{NS_stress_free} converging to a solution $\omega \in L^\infty L^2 \cap L^2 H^1$ of
\begin{equation*} 
    \begin{aligned}
        \pa_t \omega + v \cdot \na \omega - \omega \cdot \na v + \na p - \Delta \omega  & = 0, \quad\text{ in } \Omega, \\
          \div \omega & = 0, \quad \text{ in } \: \Omega \\
           \curl \omega \times n\vert_{\pa \Omega} - p n\vert_{\pa \Omega} & = 0, \quad \omega\vert_{t=0} = \omega_0  \, .
    \end{aligned}
\end{equation*}
Taking the divergence of the first equation gives $\Delta p =0$. Moreover, the normal component of the boundary condition gives $p\vert_{\pa \Omega} = 0$. Hence, $p=0$ and the above system reduces to \eqref{linear} as expected.

\section{Boundary layer analysis}  \label{secBL}
This paragraph is dedicated to the proof of Theorem \ref{thm_BL}. 

\subsection{Construction of the approximate solution} \label{subsec_approximation}
We indicate here how to construct an approximate solution of \eqref{NS}-\eqref{BC} of the form \eqref{Ansatz}, with profiles
$$ u^i = (U^i, V^i)(t,x,y), \quad u^{i,bl} = (U^{i,bl}, V^{i,bl})(t,x,z) \, . $$ 
More precisely, we shall build approximations in the form:  
\begin{equation} \label{Ansatz2}
\begin{aligned}
    u^{app}(t,x,y) & = \sum_{i=0}^N \sqrt{\nu}^i u^i(t,x,y) + \sum_{i=0}^{N-1} \sqrt{\nu}^i u^{i,bl}(t,x,y/\sqrt{\nu}) + \sqrt{\nu}^N\big(0, V^{N,bl}(t,x,y/\sqrt{\nu})\big) \\
    p^{app}(t,x,y) & = \sum_{i=0}^N \sqrt{\nu}^i p^i(t,x,y) + \sum_{i=0}^N \sqrt{\nu}^i p^{i,bl}(t,x,y/\sqrt{\nu})
\end{aligned}
\end{equation}
where the last sums on the r.h.s. start {\it a priori} from $i = 0$. We shall recover that the first two profiles are zero, hence a boundary layer of amplitude $\nu$. Note that we want $u^{app}\vert_{t=0} = u^0$, so that we impose 
$$ u^0\vert_{t=0} = u_0, \quad u^i\vert_{t=0} = 0 \quad \forall i \ge 1, \quad u^{i,bl}\vert_{t=0} = 0 \quad \forall i \ge 0. $$
Note also that we want the boundary layer correctors $u^{i,bl},p^{i,bl}$ to be localized near the boundary, in particular we want 
\begin{equation} \label{BL_infinity}
u^{i,bl},p^{i,bl} \rightarrow 0, \quad z \rightarrow +\infty      \, .
\end{equation}
To lighten notations, we also define $u^k, p^k, u^{k,bl}, p^{k,bl} \equiv 0$ for $k < 0$. We first plug the formal expansion $u^{app}$ in \eqref{BC}, and gather terms with the same powers of $\sqrt{\nu}$. The resulting equations can be put in the abstract form 
$$ \frac{1}{\nu} B^{-2}(t,x) + \frac{1}{\sqrt{\nu}} B^{-1}(t,x) + \dots \  = 0 $$
where the $B^i$ involve the $u^k$ and the $u^{k,bl}$. Imposing $B^i = 0$ for all $i$, we find 
\begin{align}
\label{BCprofiles1}  V^i\vert_{y=0} + V^{i,bl}\vert_{z=0} & = 0 \\
\label{BCprofiles2} \pa^2_z U^{bl,i}\vert_{z=0} + \pa^2_x U^{bl,i-2}\vert_{z=0} + \Delta U^{i-2}\vert_{y=0} & = 0   \,. 
\end{align}
We subsequently substitute the formal expansions $u^{app}, p^{app}$ into equation \eqref{NS}, and gather terms with the same powers of $\sqrt{\nu}$. We obtain equations of the form 
$$ \Big(  E^0(t,x,y) + \sqrt{\nu} E^1(t,x,y) + \dots \Big) + \Big( \frac{1}{\sqrt{\nu}} E^{-1,bl}(t,x,y,z) + E^{0,bl}(t,x,z)  + \dots \Big)  = 0 \, ,$$
where $z = y/\sqrt{\nu}$. Here, the $E^i$ only involves the profiles $u^k,p^k$, while the $E^{i,bl}$ involve either the profiles $u^{k,bl}, p^{k,bl}$ or mixed terms with the $u^k$ and the $u^{k,bl}$, coming from the quadratic nonlinearity. 
Sending $z$ to infinity in these equations, due to \eqref{BL_infinity}, one gets 
 $E^0 + \sqrt{\nu} E^1 + \dots = 0$ resulting in a collection of equations $E^i = 0$ involving only the $u^k$.
  We easily find that 
\begin{itemize}
    \item $u^0$ solves the Euler equation.  
    \item $u^i$, $i \ge 1$ solves in $\Omega$: 
\begin{equation} \label{eq_ui}
\pa_t u^i + u^0 \cdot \na u^i + u^i \cdot \na u^0 + \na p^i = f^i, \quad \div u^i = 0,
\end{equation}    
where $f^i$ depends on the $u^k$'s, $k \le i-1$. 
\end{itemize}
 
 Regarding the $E^{i,bl}$, we can further Taylor expand 
 $$u^k(t,x,y) = u^k(t,x,\sqrt{\nu}z) =  u^k(t,x,0) + \sqrt{\nu}z \pa_y u^k(t,x,0) + \dots $$
in the mixed terms to end up with  $\frac{1}{\sqrt{\nu}} \mathcal{E}^{-1,bl}(t,x,z) + \mathcal{E}^{0,bl}(t,x,z)  + \dots = 0$ that yields another collection of equations of the variable $z$. For instance, looking at  the $O(1/\sqrt{\nu})$ terms 
in the momentum equation, we find that $\pa_z p^{0,bl} = 0$, which combined  with \eqref{BL_infinity} yields $p^{0,bl} = 0$. Similarly, we find from the divergence-free condition that $\pa_z V^{0,bl} = 0$, hence $V^{0,bl} = 0$. We deduce from this and from \eqref{BCprofiles1} that $v_0\vert_{y=0} = 0$: this non-penetration condition completes the Euler equation and allows to determine $u^0, p^0$.  We then look at $O(1)$ terms. We find from the momentum equation that $\pa_z p^{1,bl} = 0$, hence $p^{1,bl} = 0$. From the divergence-free condition, $\pa_x U^{0,bl} + \pa_z V^{1,bl} = 0$ so that using \eqref{BL_infinity} and \eqref{BCprofiles1}
$$ V^{1,bl} = \int_z^{+\infty} \pa_x U^{0,bl}, \quad V^1\vert_{z=0} = - \int_0^{+\infty} \pa_x U^{0,bl}.  $$
The horizontal component of the momentum equation then yields 
\begin{align*} 
& \pa_t U^{0,bl} + (U^{0,bl} + U^0\vert_{y=0}) \pa_x U^{0,bl} + U^{0,bl} \pa_x U^0\vert_{y=0} \\
& + (V^{1,bl} + V^1\vert_{y=0} + z\pa_y V^0\vert_{y=0}) \pa_z U^{0,bl} - \pa_z^2 U^{0,bl} = 0 \, .
\end{align*}
It is completed with boundary condition $\pa^2_z U^{0,bl}\vert_{z=0} = 0$, see \eqref{BCprofiles2}. 
Introducing $\Omega^{0,bl} = \pa_z U^{0,bl}$, we get by differentiating the equation with respect to $z$
$$ \pa_t \Omega^{0,bl} + (U^{0,bl} + U^0\vert_{y=0}) \pa_x\Omega^{0,bl} +  (V^{1,bl} + V^1\vert_{y=0} + z \pa_y V^0\vert_{y=0}) \pa_z \Omega^{0,bl} = 0   $$
together with the homogeneous Neumann condition $\pa_z \Omega^{0,bl}\vert_{z=0}=0$. A simple estimate on this advection diffusion equation yields $\Omega^{0,bl} = 0$, and finally $U^{0,bl} = 0$. This yields $V^{1,bl} = 0$, and $V^1\vert_{y=0}$: again, this non-penetration condition completes the linearized Euler equation satisfied by $u^1, p^1$. 

Recursively, assume $i \ge 1$, and assume that we have determined 
$$X^k = \big(p^{k,bl}, U^{k-1,bl},V^{k,bl}, u^k, p^k\big), \quad k \le i$$
and that these fields are smooth and decaying enough. Then, from the momentum equation (terms of order $\sqrt{\nu}^i$), we deduce $\pa_z p^{i+1,bl} = F^i$, with $F^i$ depending on $X^k$, $k \le i$. Together with condition \eqref{BL_infinity}, this determines $p^{i+1,bl}$. From the divergence-free equation, we then get 
$\pa_x U^{i,bl} + \pa_z V^{i+1,bl} = 0$ which, together with \eqref{BL_infinity} and \eqref{BCprofiles1}, allows to express  $V^{i+1,bl}$ and $V^{i+1}\vert_{y=0}$ in terms of $U^{i,bl}$.  Writing the equation for $U^{i,bl}$ and differentiating with respect to $z$, we obtain an equation for $\Omega^{i,bl} = \pa_z U^{i,bl}$
\begin{equation}
    \pa_t \Omega^{i,bl} + U^0\vert_{y=0} \pa_x \Omega^{i,bl} + z \pa_y V^0\vert_{y=0} \pa_z \Omega^{i,bl} - \pa^2_z \Omega^{i,bl} = G^i \, , 
\end{equation}
with $G^i$ depending on $X^k$, $k \le i$. From \eqref{BCprofiles2}, we get a Neumann condition $\pa_z \Omega_i\vert_{z=0} = H^i$, with  $H^i$ depending on $X^k$, $k \le i$. We can solve this equation (see Section \ref{subsec_WP_BL} below) with zero initial data to get $U^{i,bl}$ and from there  $V^{i+1,bl}$ and $V^{i+1}\vert_{y=0}$. This last quantity provides the missing boundary condition to solve the linearized Euler equation satisfied by $u^i, p^i$ (see Section \ref{subsec_WP_BL} below). This completes the recursive construction of the profiles.  

It is straightforward to verify that $G^1 = 0$, $H^1 = 0$, so that $\Omega^{1,bl} = 0$ and so $U^{1,bl} =0$. We recover that the boundary layer part of the expansion starts at $i=2$, which corresponds to a boundary layer of amplitude $\nu$ as expected.

\subsection{Well-posedness of the profile equations} \label{subsec_WP_BL}
Let $\displaystyle \omega_0 \in C^\infty_c(\bar \Omega)$, $u_0 = \mathcal{BS}\omega_0$. Standard arguments provide existence and uniqueness of a solution $\displaystyle u^0 \in L^\infty_{loc}(\R_+, H^\infty(\Omega))$  to the Euler equation, satisfying $u^0\vert_{t=0} = u_0$. Then, given any $\displaystyle f \in L^\infty_{loc}(\R_+, H^\infty(\Omega))$ and any $\varphi = L^\infty_{loc}(\R_+,H^\infty(\pa \Omega)))$, there exists a unique solution $\displaystyle u \in L^\infty_{loc}(\R_+, H^\infty(\Omega))$ of the linearized Euler system 
\begin{align*}
    \pa_t u + u^0 \cdot \na u + u \cdot \na u^0 + \na p = f \: \text{ in } \Omega, \quad \div u = 0 \: \text{ in } \Omega, \quad u \cdot n\vert_{\pa \Omega} = \varphi. 
\end{align*}
Finally, assuming that 
$$(1+z)^n G \in L^\infty_{loc}(\R_+, H^\infty(\R \times \R_+), \quad H \in  L^\infty_{loc}(\R_+,H^\infty(\R)) \, ,$$
one can find by standard weighted estimates  a unique solution $\Omega^{bl}$ of the system 
\begin{align*}
    \pa_t \Omega^{bl} + U^0\vert_{y=0} \pa_x \Omega^{bl} + z \pa_y V^0\vert_{y=0} \pa_z \Omega^{bl} - \pa^2_z \Omega^{bl} = G, \quad \pa_z \Omega^{bl}\vert_{z=0} = H \, , 
\end{align*}
with $(1+z)^n \Omega^{bl} \in L^\infty_{loc}(\R_+, H^s(\R \times \R_+)$. 

\medskip
With these remarks, one is able to follow rigorously the iterative scheme of construction of approximate solutions indicated in the previous Section \ref{subsec_approximation}. Eventually, the approximate solution in \eqref{Ansatz2} satisfies 
\begin{equation} \label{NSapp}
    \begin{aligned}
        \pa_t u^{app} + u^{app} \cdot \na u^{app} + \na p^{app} - \nu \Delta u^{app} & = R^{app} \quad \text{in} \quad \Omega, \\
        \div u^{app} & = 0  \quad \text{in}\quad \Omega, \\
    u^{app} \cdot n\vert_{\pa \Omega} & = 0 \,  \quad \Delta u^{app} \times n\vert_{\pa \Omega} = r^{app}
    \end{aligned}
\end{equation}
where the remainders satisfy for all $p \in [1,\infty]$
\begin{equation} \label{estim_remainders}
\|R^{app}\|_{W^{k,p}(\Omega)} \lesssim \sqrt{\nu}^{N-k+1/p}, \quad \|r^{app}\|_{W^{k,p}(\pa \Omega)} \lesssim \sqrt{\nu}^{N-2} \, .
\end{equation}

\subsection{Boundary layer stability}
To conclude the proof of Theorem \ref{thm_BL}, we need to get an estimate on $\upsilon = u - u^{app}$, more precisely on  $\eta = \omega - \omega^{app}$, where $\omega$ solves the NS equation in vorticity form \eqref{eqvorticity2d}-\eqref{BCvorticity2d}, starting from $\omega\vert_{t=0} = \curl u^{app}\vert_{t=0} = \omega_0$. We find 
\begin{align*}
    \pa_t \eta + u^{app} \cdot \na \eta + \upsilon\cdot \na \omega^{app} + \upsilon\cdot \na \eta - \nu \Delta \eta & = \curl R^{app}, \quad \upsilon = \mathcal{BS} \eta \quad \text{ in } \Omega  \\
    \pa_n \eta\vert_{\pa \Omega} & = r_{app}, \quad \eta\vert_{t=0} = 0 
\end{align*}
see \eqref{NSapp}-\eqref{estim_remainders}. The keypoint is that $\na \omega^{app}$ is bounded uniformly in $\nu$. More precisely, taking into account the decay at infinity of $\omega^{app}$, one has
$$  \|(1+y)\na \omega^{app} \|_{L^\infty}  \le C. $$
A standard energy estimate then yields: 
\begin{align*} 
\frac{1}{2} \frac{d}{dt} \|\eta\|^2 + \nu \|\na \eta\|^2 & \le \|(1+y)\na \omega^{app}\|_{L^\infty} \|\frac{\upsilon}{1+y}\| \, \|\eta\| + \nu \|r^{app}\|_{H^{-1/2}(\pa \Omega)}\|\eta\|_{H^{1/2}(\pa \Omega)} \\
& + \|\curl R^{app}\| \|\eta\| \\
& \le C \|\na \upsilon\| \|\eta\| + \sqrt{\nu}^N \|\eta\|_{H^1} + \sqrt{\nu}^{N-1/2} \|\eta\|
\end{align*}
where we used Hardy's inequality $\|\upsilon/(1+y)\| \le C \|\na \upsilon\|$ to go from the first to the second line. We end up with 
\begin{align*} 
\frac{1}{2} \frac{d}{dt} \|\eta\|^2 + \nu \|\na \eta\|^2 \le C \|\eta\|^2 + \nu^{N-2}
\end{align*}
Gronwall's inequality yields the estimate of Theorem \ref{thm_BL} and concludes the proof.

\section{Final comments}
The analysis of Sections \ref{secWP2d} to \ref{secBL} give a rigorous basis to  numerical observations  made in \cite{Kerswell24a} in the context of inertial waves in a rotating flow. First, and perhaps surprisingly, the boundary conditions \eqref{BC} are sensible boundary conditions for the Navier-Stokes equations. This is not straightforward at first sight as they amount to dropping the highest order term on the tangential flow at the boundary. Second, the energy dissipation in the boundary layer is much weaker with  \eqref{BC} than with standard no-slip, or even stress-free boundary conditions. This is a nice feature when  simulating large scale flows essentially governed by an inviscid dynamics. 

\begin{figure}
\centerline{\includegraphics[height=5.6cm]{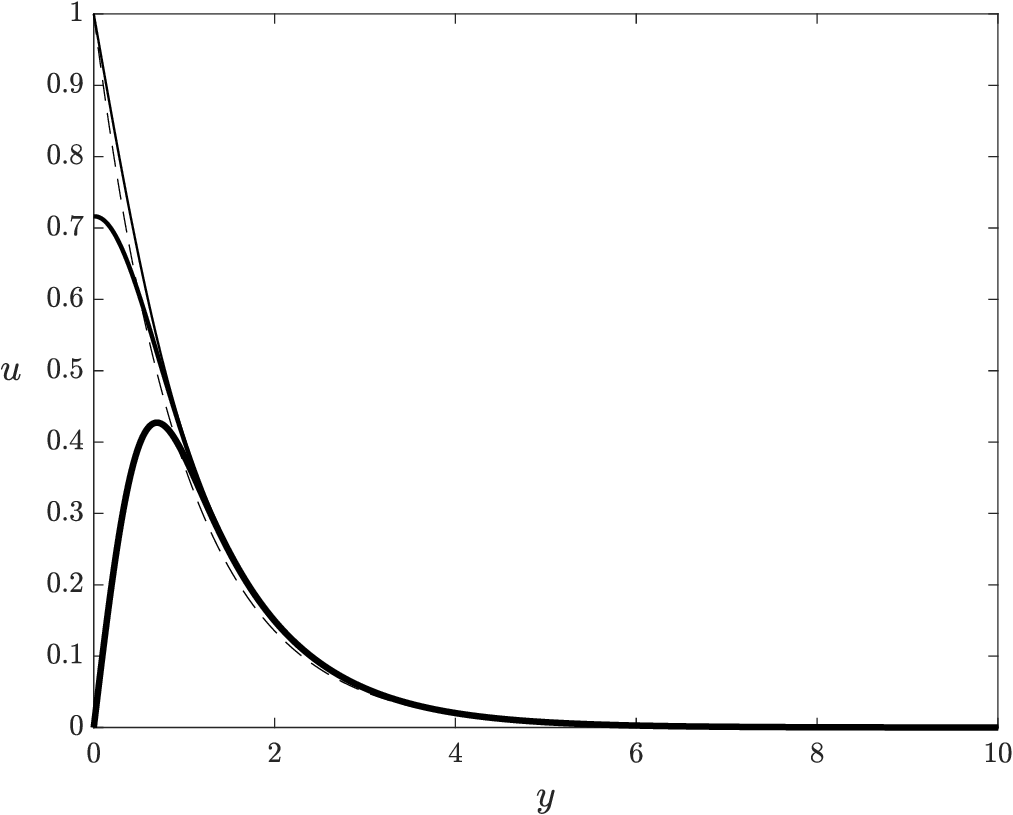}\hskip 8mm \includegraphics[height=5.6cm]{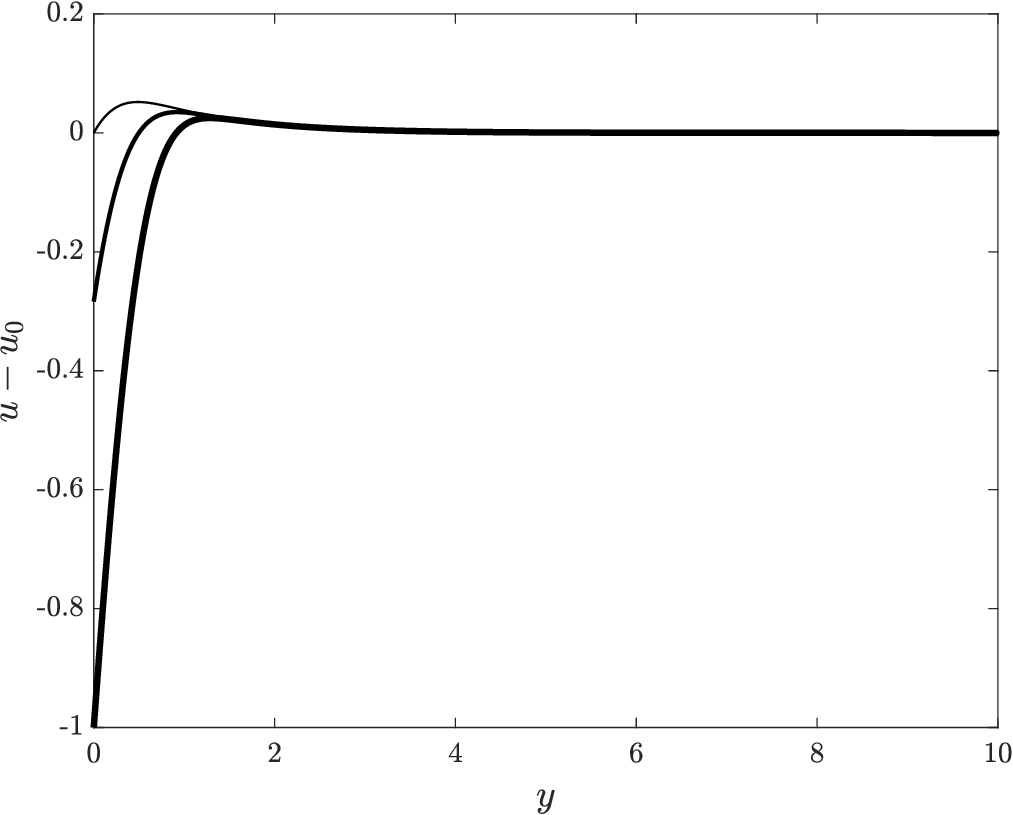}}
\caption{Numerical solution of \eqref{heat_eq}
up to time $t=1/10$ with initial condition $u_0(y)=\exp(-y) \, ,$ (in dashed line) and in decreasing line thickness, no-slip, stress-free and diffusion-free conditions. The figure on the right highlight the corrector $u-u_0\, .$
}
\label{fig1}
\end{figure}

To illustrate the different amplitudes of the boundary layers, a simple example is the case of a shear flow $(u(t,y),0)$ in a half-plane. The Navier-Stokes equation then reduces to the one dimensional heat equation 
\begin{equation}
\partial_t u = \nu \partial_{yy}u \, ,
\label{heat_eq}
\end{equation}
with boundary conditions at $y=0 \, .$
The no-slip boundary condition yields $u\vert_{y=0}=0 \, ,$ stress-free will amount to  $\partial_y u\vert_{y=0}=0 \, ,$ whereas the diffusion-free boundary condition implies $\partial_{yy} u\vert_{y=0}=0 \, .$ 
Any $(u,v)=(u_0,0)$ initial data will be a solution in the inviscid case (Euler), whereas the Navier-Stokes solution will involve corrections at the boundary of decreasing importance from no-slip to diffusion-free see figure~\ref{fig1}.

It is interesting to stress that an irrotational solution of Euler will also be a solution of Navier-Stokes (for any value of the viscosity) provided the Lions boundary condition \eqref{BC_Lions} or the diffusion-free boundary condition \eqref{BC} is used. 
Such would not be the case with the no-slip boundary condition \eqref{Noslip} or in the case of a non-flat boundary with the stress-free condition, because of the curvature term in \eqref{BC_SF3}.
For example in the case of a flow past a cylinder (or the 2D-flow past a disc), the irrotational solution of the Euler equation can be obtained using potential theory and takes the form using cylindrical $(r,\theta)$ coordinates
\begin{equation}
  u_r = U \left( 1 - \frac{a^2}{r^2} \right) \, \cos\theta \, ,\qquad
  u_\theta = -U \left( 1 + \frac{a^2}{r^2} \right) \, \sin\theta \, ,
\end{equation}
where $U$ is the flow at infinity and $a$ the cylinder radius.
This flow will be a stationary solution of the Navier-Stokes equation provided diffusion-free or Lions boundary conditions are being used. Whereas with no-slip boundary conditions, a boundary layer will form and it is well known that it will detach from the wall and emit vorticity in the downstream flow.

Another key aspect to assess the efficiency of \eqref{BC} is the way it impacts invariants of  the Euler flow. We recall that the existence of conserved quantities is related to Noether's theorem, considering the Hamiltonian nature of the  flow and its invariance by various transformation groups: translations, rotations (for further explanation, see \cite[p. 26]{Tao}). In $\R^d$, classical examples are  the kinetic energy $\frac{1}{2} \int_{\R^d} |u|^2$, the total momentum $\int_{\R^d} u$, or the total vorticity $\int_{\R^d} \omega$. In dimension $d=2$, all Casimir functionals $\int_{\R^2} f(\omega)$ are also preserved. We refer to \cite{Marchioro_Pulvirenti,Majda_Bertozzi} for further examples. In a bounded domain, under the impermeability condition $u \cdot n \vert_{\pa \Omega} = 0$, some of these quantities are still conserved as long as  symmetries of the domain are compatible with the transformation groups.

As boundary condition \eqref{BC} has a natural reformulation through vorticity, see equation \eqref{BCvorticity2d} or equation \eqref{BCvorticity3d}, it complements nicely  quantities that involve 
$\omega$. For instance, the total vorticity $\int_\Omega \omega$ is easily seen to be conserved, which is not true in the no-slip or stress-free cases. In dimension $d=2$, as seen in Paragraph \ref{nonconnected}, the circulation around obstacles is preserved. Furthermore, the Casimir functionals $\int_\Omega f(\omega)$  are under control at least for convex $f$, see the inequality  \eqref{estim_Casimir}. Finally, under \eqref{BC},  a curl-free flow  remains curl-free for all time. 

A surprising feature of equations \eqref{NS}-\eqref{BC}, however, is that the kinetic energy of the fluid may experience transient growth. Indeed, time variation of the kinetic energy is given by
$$ \frac12 \frac{d}{dt} \|u(t)\|^2 = \nu  \int_{\Omega} \Delta u \cdot u = - 2 \nu \int_\Omega |D(u)|^2 + 2 \nu \int_{\pa \Omega} D(u)n \cdot u $$
While the last integral vanishes under no-slip or stress-free boundary conditions, this is no longer the case under \eqref{BC}. To be more specific, consider the case where $\Omega = \T\times (0,1)$, corresponding to periodic boundary conditions in $x$, with walls at $z=0,1$. It can be verified that the field
$$ u_0 = (U(z),0,0), \quad U(z) = z^5 - \frac53 z^4 - \frac53 z +1.  $$
satisfies \eqref{BC2D}, which yields in this case $U''(1) = U''(0) = 0$. Moreover, it satisfies the zero flux condition : $\int_0^1 U dz = 0$. A simple calculation however shows that  
\begin{equation*}
 \nu \int_\Omega \Delta u_0 \cdot u_0 =  \frac{20 \nu}{63} > 0 \, ,  
\end{equation*} 
so that the kinetic energy of the Navier-Stokes solution with data $u_0$ will increase near the initial time. This is related to the fact that under \eqref{BC} the fluid is not isolated: stress can be provided by the solid boundary to the fluid. Mathematically, this lack of {\it a priori} $L^2$ estimate for $u$ is the reason why existence of weak Leray solutions is unclear, see Remark \ref{rem_Leray}. Nevertheless, at the level of regularity mentioned in theorems \ref{thmWP2d} and \ref{thmWP3d}, the control of the kinetic energy follows from  the control of vorticity, and at least in dimension $d=2$, the long time  behaviour can be specified, see Remark \ref{rem_exp_decay}. In the limit of vanishing viscosity, all invariants of the Euler equations are recovered.

\appendix

\section{Regularity estimates for Stokes - Proof of Lemma \ref{lemme_Stokes}} \label{appendixA}
The proof borrows ideas to  the usual proofs of elliptic regularity for the Stokes equation, when endowed with more standard boundary conditions. We notably take inspiration from \cite[chapter IV.7]{Boyer}. Here, the extra difficulty is the anisotropy of our Stokes operator, resulting in a $O(1)$ control for $\mP w$ and a $O(1/\eps)$ control for $w$. 

The keypoint is to establish \eqref{stokes_regularity_estimate} for $k=0$.  First, through introduction of a good lift $W$ such that $\div W = g$ and consideration of $w - W$, one can always assume that $g=0$. After this simplification, the variational formulation of the system becomes
\begin{equation} \label{L2stokes}
 \int_\Omega \curl w \cdot \curl \varphi - \int_{\Omega} q \div \phi + \eps  \int_\Omega \na w : \na \varphi = \int_{\pa \Omega} f_b \varphi
\end{equation}
valid for any vector field $\varphi$. With $\varphi=w$, we get
$$\|\curl w\| + \eps \|\na w\| \le C (\|f\|_{L^2} + \|f_b\|_{H^{\frac12}}  + \|w\| ) \, . $$
We can always assume that $\int_\Omega q = 0$. Under this assumption, it is well-known that there exists $\varphi_q$ satisfying 
$$ \div \varphi_q = q, \quad \varphi_q\vert_{\pa \Omega} = 0, \quad \|\varphi_q\|_{H^1} \le C \|q\|. $$
Taking $\varphi = \varphi_q$, we find after a few manipulations, using \eqref{L2stokes}
$$ \|q\| \le C (\|f\|_{L^2} + \|f_b\|_{H^{\frac12}}  + \|w\| ) \, . $$
The next step is to control $T w$ in $H^1$, for any tangential derivative $T$. We take $\varphi = T^* T \mP w$, resulting in (note that $\curl \omega = \curl \mP \omega$)
\begin{align*}
    \|\curl T \mP w \|^2  & = \int_{\Omega} \curl \mP w \cdot  [T^* , \curl] T \mP w +  \int_{\Omega} [\curl , T] \mP w \cdot  \curl  T \mP w \\
    & -\eps \int_\Omega \na w \cdot [\na ,  T^*] T \mP w + \eps \int_{\Omega} T \na w \cdot \na T \mP \omega \\
    & +\int_\Omega q \,  [\div , T^*] T \mP w  +  \int_\Omega T q \,  [\div , T] \mP w  +    \int_{\pa \Omega} T f_b \cdot  T \mP w  \, .
\end{align*}
Note that the last term at the right-hand side obeys the bound 
$$\left|\int_{\pa \Omega} T f_b \cdot  T \mP w\right| \le C \|T f_b\|_{H^{-1/2}(\pa \Omega)} \|T \mP w\|_{H^{1/2}(\pa \Omega)} \le C' \|f_b\|_{H^{1/2}(\pa \Omega)} \|T \mP w\|_{H^1(\Omega)} \, .$$ 
The other integrals are treated classicaly, and we get  
\begin{align*}
    \|\curl T \mP w \|^2  & \le  C \Big(   \|w\|_{H^1} \|T \mP w\|_{H^1}  + \eps \|Tw\|_{H^1} \|T \mP w\|_{H^1}  \\
    & +  \|q\| \|T\mP w\|_{H^1} + \|Tq\| \|w\|_{H^1} + \|f_b\|_{H^{1/2}(\pa \Omega)} \| \mP w \|_{H^1} \Big) \, .
\end{align*}
As $\mP w \cdot n\vert_{\pa \Omega} = 0$, we find that $\|T\mP\omega \cdot n\|_{H^{1/2}(\pa \Omega)} = 
\|\mP\omega \cdot Tn\|_{H^{1/2}(\pa \Omega)} \le C \|\omega\|_{H^1}$, and so 
$$ \|\na T \mP \omega\| \le C (\|\curl T \mP \omega\| + \|T\mP\omega \cdot n\|_{H^{1/2}(\pa \Omega)}) \le C' (\|\curl T \mP \omega\| + \|\omega\|_{H^1})  \, .$$
Combining this inequality with the previous one, playing with Young's inequality, we infer that 
\begin{equation} \label{estimTpw}
    \|T \mP w \|^2_{H^1} \le C \big(  \eps^2 \|Tw\|_{H^1}^2 +  \|w\|_{H^1}^2 +  \|Tp\|^2_{L^2} +  \|f\|^2 + \|f_b\|^2_{H^{\frac12}(\pa \Omega)} \big) \, .
\end{equation}
For the next estimate, we take $\varphi = \eps T^* T w$: 
\begin{align*}
\eps^2 \|\na T w\|^2 & = -\eps \int_{\Omega} \curl \mP w \cdot [\curl , T^*] T w  -\eps \int_{\Omega} T\curl \mP w \cdot \curl T w \\
& + \eps^2 \int_\Omega \na w \cdot [T^* , \na] T  w + \eps^2 \int_{\Omega} [\na , T] w \cdot \na T \omega \\
    & + \eps \int_\Omega q \,  [\div , T^*] T  w  +  \eps \int_\Omega T q \,  [\div , T]  w  +  \eps   \int_{\pa \Omega} T f_b \cdot  T  w  \, .
\end{align*}
The only delicate term at the right-hand side is the second one, for which it is crucial to notice that 
\begin{align*}
\curl T w = \curl T \mP w + \curl T \mQ w = \curl T \mP w  + [\curl, T] \mQ \omega 
\end{align*}
so that 
$$ \|\curl T w\| \le C ( \|T \mP w\|_{H^1} + \|\mQ \omega\|_{H^1} ) \le  C' ( \|T \mP w\|_{H^1} + \|\omega\|_{H^1}) \, . $$
We thus get
\begin{equation} \label{estimTw}
\eps^2 \|T w\|_{H^1}^2  \le    C \big( \eps \|T \mP w\|_{H^1}^2 +  \|w\|_{H^1}^2 + \eps^2 \|Tq\|^2 +  \|f\|^2 + \|f_b\|^2_{H^{\frac12}(\pa \Omega)} \big)  
\end{equation}
The last tangential estimate  to get is the one for $Tq$. We therefore introduce $\Phi_q$ such that 
$$ \Delta \Phi_q = T q, \quad  \Phi_q\vert_{\pa \Omega} = 0$$
We have in particular (notice that $\int_\Omega Tq \, \Phi_q = - \int_{\Omega} q \,  T^* \Phi_q$)
$$\|\Phi_q\|_{H^1} \le C \|q\|, \quad \|\Phi_q\|_{H^2} \le C \|T q\|\, .$$ 
We then  take $\varphi = T^* \na \Phi_q$. After manipulations similar to the previous ones, we find 
\begin{equation} \label{estimTp}
    \|T q\| \le C \big( \eps\|T w \|_{H^1} + \|q\| + \|w\|_{H^1} + \|f\| + \|f_b\|_{H^{\frac12}(\pa \Omega)}\big) \, .
\end{equation}
Gathering of \eqref{estimTpw},\eqref{estimTw},\eqref{estimTp} yields
$$ \|T\mP w\|_{H^1} + \eps \|w\|_{H^1} + \|Tq\| \le C \big(\|f\| + \|f_b\|_{H^{\frac12}(\pa \Omega)} + \|w\|_{H^1}\big) \, .$$
The final step is to obtain a bound on the normal derivatives $Nq$ and $N N w$ in $L^2$. Using the divergence-free condition, $Nw \cdot n$ can be expressed in terms of tangential derivatives of $w$ (and lower order terms), hence it is bounded in $H^1$. From there, we can consider the normal component of the Stokes equation in \eqref{stokes_eq}, and express (up to lower order terms) $Np$ with  derivatives of $Tw$, derivatives of $Nw \cdot n$ and $f$. Hence, it is bounded in $L^2$. Eventually, considering the tangential components of the Stokes equation yields a similar $H^1$ control of the tangential part of $Nw$.

 Once such $H^2$ regularity estimate is obtained, the $H^{k+2}$ estimates, $k \ge 1$, can be obtained inductively. Typically, to go from $H^2$ to $H^3$ regularity , one proceeds as follows: 
\begin{itemize}
\item Let $T$ an arbitrary tangential derivative. Applying $T$ to system \eqref{stokes_eq}, one sees that the field $Tw$ still obeys a Stokes system of type \eqref{stokes_eq}:  $Tf$ and $Tf_b$ replace $f$ and $f_b$, and there is an extra source term, resp. divergence term, resp.   boundary term,  coming from commutators. The $H^2$ estimate shows that these extra terms belong to $L^2$, resp. $H^1$, resp. $H^{\frac12}$. One can then apply the $H^2$ regularity estimate to control $Tw$  in $H^2$ and $Tp$ in $H^1$. 
\item Regarding the remaining normal derivatives $NNp$  and $NNN w$, one  can apply $N$ to \eqref{stokes_regularity_estimate} and deduce a similar Stokes type equation on $Nw$. But one can no longer differentiate the boundary condition, as $N$ is normal to the boundary. Instead, the trick is to deduce from the equalities
$$ (\curl w \times n + \eps \pa_n w - q n)\vert_{\pa \Omega} = f_b, \quad \div w\vert_{\pa \Omega} = g\vert_{\pa \Omega}  \, ,$$
that $Nw$ satisfies an inhomogeneous Dirichlet condition, involving tangential derivatives $Tw$. By the first step, this Dirichlet datum belongs to $H^{\frac32}(\pa \Omega)$. By the classical $H^2$ estimate for the Stokes equation with Dirichlet condition, this finally implies that $Nw \in H^2$, $Np \in H^1$. 
\end{itemize}

\bigskip
\bigskip
\noindent
{\em Statements and declarations.}

\medskip 
\noindent
Data availibility statement: The authors declare that the data supporting the findings of this study are available within the paper. 

\medskip
\noindent 
The authors have no competing interests to declare that are relevant to the content of this article.
They did not receive support from any organization for the submitted work.

\bibliographystyle{siam}
\bibliography{papier.bib}

\end{document}